\documentclass{nmeauth}
\usepackage{}

\usepackage{epsfig}
\usepackage{subcaption} 
\usepackage[abs]{overpic}
\usepackage{bm}

\makeatletter

\newcommand{\im}{\mathrm{i}}
\newcommand{\diff}{\mathrm{d}}
\newcommand{\rank}{\mathrm{rank}}

\newcommand{\spanset}{\mathrm{span}}

\newcommand{\tol}{\mathrm{tol}}
\newcommand{\gap}{\mathrm{gap}}

\newcommand{\ssrrs}{RSRR}
\newcommand{\ssrrm}{SSRR}

\newtheorem{thm}{Theorem}

\theoremstyle{definition}
\newtheorem{example}{Example}

\newcommand\BibTeX{{\rmfamily B\kern-.05em \textsc{i\kern-.025em b}\kern-.08em
T\kern-.1667em\lower.7ex\hbox{E}\kern-.125emX}}

\definecolor{shadecolor}{rgb}{0.9,0.9,0.9}

\makeatother

\def\volumeyear{2016}

\begin{document}

\runningheads{J.~Xiao, \textit{et al}}{Rational interpolation approach and resolvent sampling scheme}

\title{Solving large-scale nonlinear eigenvalue problems by rational interpolation approach and resolvent sampling based Rayleigh-Ritz method}

\author{Jinyou Xiao\affil{1}\comma\corrauth, Chuanzeng Zhang\affil{2}, Tsung-Ming Huang\affil{3}, Tetsuya Sakurai\affil{4}}

\address{\affilnum{1}School of Astronautics, Northwestern Polytechnical University, Xi'an 710072, China\break
\affilnum{2}Department of Civil Engineering, University of Siegen, D-57068 Siegen, Germany\break
\affilnum{3}Department of Mathematics, National Taiwan Normal University, Taipei 116, Taiwan\break
\affilnum{4}Department of Computer Science, University of Tsukuba, Tsukuba, 305–8573, Japan
}

 \corraddr{School of Astronautics, Northwestern Polytechnical University, Xi'an 710072, China.
Email: xiaojy@nwpu.edu.cn
}

\begin{abstract}

Numerical solution of nonlinear eigenvalue problems (NEPs) is frequently encountered in computational science and engineering. The applicability of most existing methods is limited by matrix structures, property of eigen-solutions, size of the problem, etc. This paper aims to break those limitations and to develop robust and universal NEP solvers for large-scale engineering applications. The novelty lies in two aspects. First, a rational interpolation approach (RIA) is proposed based on the Keldysh theorem for holomorphic matrix functions. Comparing with the existing contour integral approach (CIA), the RIA provides the possibility to select sampling points in more general regions and has advantages in improving accuracy and reducing computational cost. Second, a resolvent sampling scheme using the RIA is proposed for constructing reliable search spaces for the Rayleigh-Ritz procedure, based on which a robust eigen-solver, denoted by \ssrrs{}, is developed for solving general NEPs. \ssrrs{} can be easily implemented and parallelized. The advantages of the RIA and the performance of \ssrrs{} are demonstrated by
a variety of benchmark and practical problems.

\end{abstract}

\keywords{Eigenvalue problems; nonlinear solvers; Rayleigh-Ritz procedure; finite element methods; boundary element methods}

\maketitle

\section{Introduction}
\label{s-intro}
\vspace{-2pt}

Robust numerical solution of large-scale nonlinear eigenvalue problems (NEPs) is of great importance in many fields of computational science and engineering.
However, besides the quadratic eigenvalue problem that has been well-studied in the dynamic analysis \cite{TM01}, the numerical solution of more general NEPs is so far a challenging task \cite{MS11, effenberger2013robust, van2015rational}.

This work aims to develop efficient numerical methods for solving large-scale NEPs of the general form
\begin{equation}\label{eq-gnep}
T(\lambda)v=0,
\end{equation}
where $T(z) \in \mathbb{C}^{n \times n}$ is a matrix-valued function depending on a parameter $z \in \mathcal{D} \subset \mathbb{C}$.
In particular, we consider the NEPs in the \emph{finite element method} (FEM) and the \emph{boundary element method} (BEM), as these two methods are extensively used in real applications, either independently or in coupled manners \cite{OBZ15,Steinbach14K}, and NEPs of them reflect the main bottlenecks in the current development of numerical methods for NEPs. However, the methods developed in this paper are rather general and by no means limited to the NEPs in the FEM and BEM.

Generally speaking, although there exist a number of numerical methods for NEPs in literature, most of them are restricted by matrix structures, properties of eigenvalues, computational costs, etc. Methods that can robustly and reliably calculate all the eigenpairs within a given region, and at the same time, can be applied to large-scale solutions are still lacking \cite{abdel1994safeguarded, dumont2007solution, effenberger2013robust, EC13, salas2014spectral, van2015rational, lu2015pade}.

In recent years, the \emph{contour integral approach (CIA)} based on the contour integrals of the probed resolvent $T(z)^{-1}U$ has attracted much attention \cite{sakurai2003projection, polizzi2009density, AS09, Beyn12, YS13, gavin2013non, austin2015computing, stefan2015zolotarev}; here, the constant matrix $U \in \mathbb{C}^{n \times L}$ consists of $L$ linear independent column vectors that are used to probe $T(z)^{-1}$.
The CIA is first developed for solving generalized eigenvalue problems \cite{sakurai2003projection, polizzi2009density}, and later, extended to solving NEPs in \cite{SAT09} and \cite{Beyn12} using the Smith form and Keldysh's theorem for analytic matrix-valued functions, respectively. Since \cite{SAT09} and \cite{Beyn12} result in similar algorithms, we consider the block Sakurai-Sugiura (SS) algorithm in \cite{SAT09}. This algorithm transforms the original NEP into a small-sized generalized eigenvalue problem involving block Hankel matrices by using the monomial moments of the probed resolvent $T(z)^{-1} U$ on a given contour.
It becomes unstable and inaccurate when higher order contour moments of $T(z)^{-1} U$ are used.

To circumvent the problem, the \underline{S}akurai-\underline{S}ugiura method with the \underline{R}ayleigh-\underline{R}itz projection (\ssrrm{}) has been recently proposed \cite{YS13}.
In \ssrrm{} the eigenspaces are constructed by a matrix $M \in \mathbb{C}^{n \times K\cdot L}$ collecting all the monomial moments of $T(z)^{-1} U$ up to a given order $K$. We refer to this as the \emph{resolvent moment scheme} (or simply, \emph{moment scheme}) for generating eigenspaces.
Numerical results in \cite{YS13} show that the \ssrrm{} algorithm has much better stability and accuracy than the SS algorithm when a small value of $L$ is used. Besides, the \ssrrm{} algorithm inherits all the merits of the CIA: (1) it can simultaneously compute all the eigenvalues (and the associated eigenvectors) within a simple closed contour; and (2) the most computationally intensive part, the computation of $T(z)^{-1} U$ at a series of sampling points $z_i, i=0,1,\dots,N-1$, is easily parallelizable.

However, the \ssrrm{} algorithm is likely to fail when high order moments of $T(z)^{-1} U$ are used to construct search spaces, because these moments tend to be linearly dependent with the increase of order $K$, and finally it is unlikely to obtain proper eigenspaces no matter how large $K$ is!
Instead of using a large $K$, one can ensure the quality of the constructed eigenspaces by using a large $L$; the FEAST algorithm, due to Polizzi \cite{polizzi2009density}, is such an example with $K=0$. While bear in mind that the main computational cost of the \ssrrm{} algorithm increases with $L$, thus a small value of $L$ is always preferable in solving large-scale problems.

Motivated by the attractive advantages of the CIA, the goal of this paper is to generalize the CIA and develop more robust and efficient solvers for large-scale engineering NEPs. There are two main contributions. First,
We propose a \emph{rational interpolation approach} (RIA) for solving NEPs. It provides a more general framework and includes the CIA as a special case where the sampling points $z_i\, (i=0,1,\dots,N-1)$ are placed on a contour enclosing the domain of interpolation. The SS and \ssrrm{} algorithms based on the CIA are generalized to the framework of the RIA. This extension offers more freedom in the selection of the sampling points $z_i$, which can be further exploited to improve the accuracy and computational efficiency.

Then, based on the RIA, we propose a \emph{resolvent sampling scheme} (or simply, \emph{sampling scheme}) for constructing approximate eigenspaces and a \underline{R}esolvent \underline{S}ampling based \underline{R}ayleigh-\underline{R}itz method, abbreviated as \ssrrs{}, for solving general NEPs. Compared with the moment scheme, the sampling scheme has three advantages: (1) it does not use any moment of $T(z)^{-1}U$ and thus effectively circumvents the possible failure of the \ssrrm{} algorithm; (2) it generates a larger and more reliable subspace; and (3) it allows to put the sampling points close to the eigenvalues so as to further improve the accuracy of the eigen-solutions. These advantages make the \ssrrs{} algorithm much more robust and accurate than the \ssrrm{} algorithm.

Below is a brief review of the other state-of-the-art methods for solving general NEPs.
Essentially, two types of methods are in use \cite{VRK13}: methods that directly deal with the NEPs, and methods that transform the NEPs into polynomial or rational eigenvalue problems and then solve via linearization.
For the former type of methods, we find examples as the residual inverse iteration method \cite{Neu85}, the Jacobi-–Davidson method \cite{Voss07}, and the block Newton method \cite{kressner2009block}.
The Jacobi-–Davidson method is perhaps the most promising one, but there are still a few issues having crucial influence on the performance of the method in practical applications.
For example, how to effectively accelerate the solution of the correction equation and control its accuracy \cite{hochstenbach2009controlling}, and how to inhibit the method from repeatedly converging towards the same eigenvalues \cite{EC13}.

Linearization is a standard approach in solving quadratic and rational eigenvalue problems in the FEM \cite{TM01, effenberger2012linearization}. Its application to the NEPs in the BEM can be found in \cite{kirkup1993solution, ARY95}. Recently, linearization by the polynomial or rational approximation of $T(z)$ has attracted increasing attentions; see, e.g., \cite{EK12, van2015rational}. In particular, the compact rational Krylov method in \cite{van2015compact} exploits the structure of the linearization pencils by using a generalization of the compact Arnoldi decomposition. As a result, the extra memory and computational costs of the linearization can be greatly reduced for large-scale problems.
However, for this class methods there is still a need to construct linearizations that reflect the
structure of the given matrix polynomial and to improve the stability of the linearizations \cite{mackey2006structured}.

Finally, we notice that previous nonlinear eigensolvers would be problematic in dealing with large NEPs arising in the BEM and the related methods \cite{EK12, EKSU12}, due to the fact that the matrix $T(z)$ is typically complex, dense and unstructured, and the evaluation of $T$ itself and the operations with $T$ (e.g., applying to vectors, solving linear systems, etc) are often computationally expensive.
For example, when the Jacobi-–Davidson method is used, the computational costs for the repeated formulation and solution of the correction equations during the iteration process would quickly become unaffordable. When using the linearization methods, a key problem is how to store all the dense coefficient matrices of $T(z)$ in, for example, a polynomial basis; to our knowledge, even using the current fast BEM techniques would require huge memory for large-scale problems.
On the contrary, the \ssrrs{} algorithm only involves the solution of a series of mutually independent equations $T(z_i)^{-1} U$, thus each matrix $T(z_i)$ is computed and used only once.

The rest of this paper is organized as follows. In Section \ref{S-basics}, some typical NEPs in the FEM and BEM, as well as the Keldysh theorem for the NEPs are briefly reviewed for completeness.
In Section \ref{S-RI} the RIA and a SS-type algorithm SS-RI are developed.
In Section \ref{S-Rayleigh}, the sampling scheme is first developed based on the RIA, and then the main eigensolver, \ssrrs{}, is presented.
In Section \ref{S-nin} the theoretical benefits of the RIA and the performance of the algorithms SS-RI and \ssrrs{}
are numerically studied.  
In Section \ref{S-ne-realapps}, the performance of \ssrrs{} is further demonstrated by two large-scale practical problems in FEM and BEM, respectively.
The essential conclusions of the paper is presented in Section \ref{S-conclusions}.

\section{Basics} \label{S-basics}

Some representative NEPs in the FEM and BEM are first briefly summarized, which will be used as target problems in the numerical examples. Then, the Keldysh theorem for the NEPs is reviewed.

\subsection{Typical NEPs}
\label{S-S-typicalNEPs}

The most widely studied NEP is the \emph{quadratic eigenvalue problem} in the dynamic analysis of structures, which takes the form
\begin{equation}\label{eq-QEP}
T(\lambda) = \lambda^2 M + \lambda C + K_{\rm s}.
\end{equation}
Typically, the stiffness matrix $K_{\rm s}$ and the mass matrix $M$ are real symmetric and positive
(semi-)definite, and the damping matrix $C$ is general. However, for FE models of rotating
machinery, the matrices $K_{\rm s}$ can be nonsymmetric due to the influence of the gyroscopic and circulatory forces \cite{quraishi2014solution}.

The second type of NEP is the \emph{rational eigenvalue problem}. Examples include the case
\begin{equation}\label{eq-REP-string}
T(\lambda) = \lambda M - K_{\rm s} + \sum_{j=1}^{J}{\lambda \over \sigma_j - \lambda} C_j
\end{equation}
that occurs in the study of the free vibration of plates with elastically attached masses \cite{SSI06}
or vibrations of fluid solid structures \cite{conca1989existence}, and the case
\begin{equation}\label{eq-REP-damping}
T(\lambda) = \lambda^2 M + K_{\rm s} + \lambda\left(a_0 + \sum_{j=1}^J {a_j\over \lambda + b_j} \right)K_{\rm v}
\end{equation}
that arises in the modal analysis of structures with viscoelastic damping treatment (such as composite structural materials, active control and damage tolerant systems in airplane, rocket, etc) \cite{trindade2000modeling}.
In \eqref{eq-REP-string} the parameters $\sigma_j$ are given by the boundary conditions, $C_j$ are damping matrices. In \eqref{eq-REP-damping} $a_j$ and $b_j$ are the relaxation parameters of the given damping model,
$K_{\rm v}$ is the unit viscoelastic stiffness matrix. Note that the above two rational problems can be turned into \emph{polynomial eigenvalue problems} by multiplying with
an appropriate scalar polynomial in $\lambda$.

The third type of NEP arises from the FE discretization of the boundary value problems involving the Maxwell equation in electromagnetic modeling of waveguide loaded accelerator cavities \cite{liao2010nonlinear}. It takes the form
\begin{equation}\label{eq-nep-gun}
T(\lambda) = K_{\rm s}- \lambda^2 M + \im \sum_{j=1}^J \sqrt{\lambda^2-\kappa_j^2} \,C_j,
\end{equation}
with $C_j$ being the damping matrices, and $\im = \sqrt{-1}$. The nonlinearity of the NEPs \eqref{eq-REP-string} and \eqref{eq-nep-gun} are caused by the nonlinear dependence of boundary conditions on $\lambda$, which is distinct from the case \eqref{eq-REP-damping} in which the nonlinearity stems from the dependence of material properties on $\lambda$.

Lastly, we mention a more general type of NEP in which the dependence of $T$ on $\lambda$ is complicated or not explicitly known. An example is the NEP in the thermoacoustic simulations involved in the stability analysis of large combustion devices \cite{nicoud2007acoustic},
\begin{equation}\label{eq-gnep-fem}
T(\lambda) = A + \lambda B(\lambda) + \lambda^2 C - D(\lambda),
\end{equation}
where matrices $A$, $B$, $C$ and $D$ are sparse matrices obtained from the FEM discretization of the differential equations and boundary conditions. $B$ corresponds to the acoustic impedance boundaries of the combustion chamber, which becomes $\lambda$-dependent when the impedance is $\lambda$-dependent. $D$ accounts for the interaction between the sound pressure and the flame, and it is a nonlinear function of $\lambda$ whose expression is generally not available.

Another example in this aspect is the NEPs in the BEM. The semi-discretized form of the direct boundary integral equation is $H(\lambda)u(\lambda) = G(\lambda)q(\lambda)$,
where $H(\lambda)$ and $G(\lambda)$ are complex square matrices, $u$ and $q$ are vector collections of the nodal displacement and traction components; see, e.g. \cite{CWX15}.
The matrix $T(\lambda)$ of the cooresponding NEP consists of the columns of the matrices $H(\lambda)$ and $G(\lambda)$ according to the given boundary conditions. In general, the entries of $H(\lambda)$ and $G(\lambda)$ are distinct functions of $\lambda$ whose expressions can not be obtained explicitly. For instance,
in the acoustic Nystr\"om BEM in \cite{CWX15},
the entries of $H(\lambda)$ and $G(\lambda)$ associated with well-separated nodes $x_i$ and $y_j$ are given by the values of double- and single-layer integral kernels,
\begin{equation} \label{eq-nyes}
H_{ij}(\lambda) = {\partial \over \partial n_{y_j}} {\exp(\im \lambda r_{ij}) \over r_{ij}}, \quad G_{ij}(\lambda) =  {\exp(\im \lambda r_{ij}) \over r_{ij}},
\end{equation}
where $r_{ij} = |x_i-y_j|$, $n_{y_j}$ is the outward normal of the boundary at node $y_j$. When the nodes $x_i$ and $y_j$ are not well-separated, the element integrals correspond to $x_i$ and the boundary element of $y_j$ become singular or nearly singular, and the expressions of $H_{ij}(\lambda)$ and $G_{ij}(\lambda)$ in terms of $\lambda$ often can not be obtained in close forms.

\subsection{Basic theory of NEPs} \label{S-neptheory}

In all the previously mentioned NEPs, $T(\lambda)$ is a nonlinear matrix-valued function of $\lambda$. In order to build our methods on a solid mathematical foundation, we further assume that $T(\lambda)$ is holomorphic in a neighborhood of all the eigenvalues. For the NEP in \eqref{eq-REP-string}, this condition implies that the eigenvalues $\lambda$ should not be infinitely close to the parameters $\sigma_j$, which is indeed inherently satisfied in most practical problems.

Now we consider a holomorphic matrix-valued function $T(z) \in \mathbb{C}^{n \times n}$ defined in an open domain $\mathcal {D} \in \mathbb{C}$, and assume that the determinant $\det T(z)$ does not vanish identically.
We intend to search for all the eigenvalues within a compact set $C \subset \mathcal {D}$ and the associated eigenvectors. Our method is motivated by the following theorem regarding the relation between the eigenvalues and the resolvent $T(z)^{-1}$ (see \cite{Beyn12}, Corollary 2.8).
\begin{thm} \label{thm1}
Let $C \subset \mathcal {D}$ be a compact set containing a finite number $n_C$ of different eigenvalues $\lambda_k\, (k=1, \cdots, n_C)$, and let
\begin{equation} \label{eq-VC}
V_C = \left ( v^{l,k}_j, \quad 0 \le j \le \mu_{l,k} -1, \, 1 \le l \le \eta_k, \, k=1, \cdots, n_C \right)
\end{equation}
and
\begin{equation} \label{eq-WC}
W_C = \left ( w^{l,k}_j, \quad 0 \le j \le \mu_{l,k} -1, \, 1 \le l \le \eta_k, \, k=1, \cdots, n_C  \right)
\end{equation}
be the corresponding \emph{canonical systems of generalized eigenvectors} (CSGEs) of $T(z)$ and $T(z)^H$, respectively. Then there exists a neighborhood $C \subset \mathcal {U} \subset \mathcal {D}$ and a holomorphic matrix-valued function $R_C: \Omega \rightarrow \mathbb{C}^{n \times n}$ such that for all $z \in \mathcal {U} \setminus \{\lambda_1, \cdots, \lambda_{n_C }\}$, it holds
\begin{equation}\label{eq-thm1}
T(z)^{-1} = \sum^{n_{C}}_{k=1}  \sum^{\eta_k}_{l=1}  \sum^{\mu_{l,k}}_{j=1}  ( z-\lambda_k )^{-j}  \sum^{\mu_{l,k}-j}_{m=0}    v^{l,k}_m \left(w^{l,k}_{\mu_{l,k} -j-m}\right)^{H} + R_C(z).
\end{equation}
\end{thm}

In expressions \eqref{eq-VC}, \eqref{eq-WC} and \eqref{eq-thm1}, $\eta_k$ and $\mu_{l,k}$ represent the dimension of the nullspace of $T(\lambda_k)$ and the $l$th partial multiplicity of $T(z)$ at $\lambda_k$, respectively. We refer the readers to \cite{Beyn12, YS13} for more detailed description of the related notations. The superscript $H$ denotes the conjugate transpose.

The summation in \eqref{eq-thm1} can be recast into the matrix form as
\begin{equation}\label{eq-VPhiW}
\sum^{n_{C}}_{k=1}  \sum^{\eta_k}_{l=1}  \sum^{\mu_{l,k}}_{j=1}  ( z-\lambda_k )^{-j}  \sum^{\mu_{l,k}-j}_{m=0}    v^{l,k}_m \left({w}^{l,k}_{\mu_{l,k} -j-m}\right)^{H}
 = V_C  \Phi(z)  W_C^H,
\end{equation}
where the matrices $V_C$ and $W_C$ are consist of the CSGEs of $T(z)$ and $T(z)^H$, respectively; the matrix-valued function $\Phi(z)$ is given by
\begin{equation} \label{eq-MPhiz}
\begin{split}
&\Phi(z)= \begin{bmatrix}
  \Phi_1(z) &   &   \\
    & \ddots &   \\
    &   & \Phi_{n_C}(z) \\
\end{bmatrix}, \quad
\Phi_k(z) = \begin{bmatrix}
  \Phi_{k}^{1}(z) &   &   \\
    & \ddots &   \\
    &   & \Phi_{k}^{\eta_k}(z) \\
\end{bmatrix},\\
&\Phi_{k}^{l}(z)= \begin{bmatrix}
           ( z-\lambda_k )^{-1}  & ( z-\lambda_k )^{-2} & \cdots  &  ( z-\lambda_k )^{-\mu_{l,k}} \\
             & \ddots & \ddots &  \vdots \\
             &   & ( z-\lambda_k )^{-1} & ( z-\lambda_k )^{-2} \\
             &   &   & ( z-\lambda_k )^{-1}  \\
         \end{bmatrix}.
\end{split}
\end{equation}
By using \eqref{eq-VPhiW}, the resolvent $T(z)^{-1}$ can be expressed more concisely as
\begin{equation}\label{eq-thm1-mf}
T(z)^{-1} = V_C  \Phi(z)  W_C^H + R_C(z).
\end{equation}

Now we consider a domain enclosed by a Jordan cure $\mathcal {C}$ in the complex plane. We are interested in all the eigenvalues (and the associated eigenvectors) inside $\mathcal {C}$.
Denote by $n_{\mathcal {C}}$ the number of mutually different eigenvalues inside $\mathcal {C}$, and by $\bar n_{\mathcal {C}}$ the total number of eigenvalues counting the algebraic multiplicity, i.e.,
$\bar n_{\mathcal{C}} = \sum^{n_{\mathcal{C}}}_{k=1}  \sum^{\eta_k}_{l=1} \mu_{l,k}$.
The corresponding CSGEs of $T(z)$ and $T(z)^H$ are denoted by $V_\mathcal {C}$ and $W_\mathcal {C}$, respectively.

\section{Rational interpolation approach} \label{S-RI}

The RIA developed in this section is a generalization of the CIA by Sakurai and Sugiura \cite{AS09} that aims to circumvent the inherent limitations of the CIA and enhance the accuracy and computational efficiency of the resulting numerical algorithms. For further distinction,
the SS algorithm \cite{AS09} based on the CIA will be denoted by \emph{SS-CI}.

In Section \ref{S-S-motivation} the theoretical limitations of the CIA is discussed. Then, in Sections \ref{S-S-ria} and \ref{S-S-rialgorithm} the RIA and a SS-type algorithm using the RIA, denoted by \emph{SS-RI}, are described. The SS-RI algorithm is a straight-forward application of RIA. It reflects the potential advantages of the RIA, however, like SS-CI, its performance is sensitive to parameter choice and thus it is not suitable for large-scale practical applications. This shortcoming further motivates the Rayleigh-Ritz reformulation in Section \ref{S-Rayleigh}.

\subsection{Motivation} \label{S-S-motivation}

\begin{figure}[hbt]
\centering
\epsfig{figure=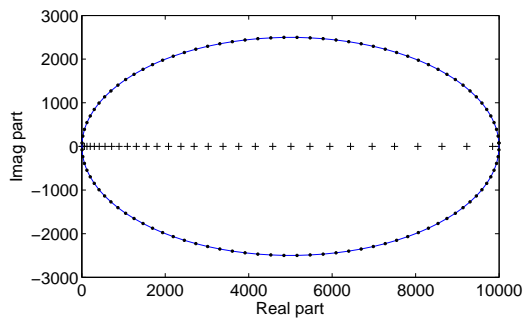,width=0.5\textwidth}
\caption{The eigenvalues ($\scriptstyle\pmb{+}$), contour (the solid line) and contour quadrature points ($\bullet $) for the SS-CI algorithm in Example \ref{ex-string} of Section \ref{S-S-nin-ss} }
\label{fig-string-evs}
\end{figure}

Figure \ref{fig-string-evs} diagrammed a situation in which the interested eigenvalues lie in a real interval. When using SS-CI, the contour should be chosen as a simple closed curve enclosing the interval. Both circles and ellipses are frequently used. Theoretical results based on rational filter theory indicate that circular contours should lead to better accuracy since in this case the filters are more closer to the indicator function of the interval \cite{stefan2015zolotarev, van2015designing}. Numerical results, however, demonstrate that flat elliptical contours often achieves better results. More interestingly, we found that the ``quadrature'' points can even lie in the interested interval. This scenario is obviously not covered by the CIA, but it works well and sometimes might be advantageous in saving computational cost.
Furthermore, for more general domains (like a rectangular domain), we found that sampling points inside the domain can also be used, and by using sampling points close to the eigenvalues, the accuracy of the eigen-solutions can often be considerably improved; see Section \ref{S-S-accuracy-gun}.

To be able to interpret the above findings, the rational interpolation approach (RIA) is developed based on Theorem \ref{thm1}.

\subsection{Theoretical derivation} \label{S-S-ria}

Theorem \ref{thm1} shows that the eigenvalues of $T(z)$ are the poles of its resolvent $T(z)^{-1}$. This inspires us to solve NEPs using the pole-finding methods. To construct practical algorithms, we consider the
following reduced matrix-valued function $F(z)$ having the same poles as $T(z)^{-1}$,
\begin{equation}\label{eq-F}
F(z) = U^H T^{-1}(z)U,
\end{equation}
where $U \in \mathbb{C}^{n \times L}$ is a constant matrix of full rank used to condense $T(z)^{-1}$. The use of $F(z)$ instead of $T(z)^{-1}$ will greatly reduce the computational time and storage of the following algorithms, which is essential for large-scale problems. For small-scale problems $U$ can be discarded, i.e., directly letting $F(z) = T^{-1}(z)$.
To properly extract all the poles of $T(z)^{-1}$ with right multiplicities, the number $L$ has to be at least equal to the maximal algebraic multiplicity of eigenvalues of $T(z)$ in $\mathcal {C}$ \cite{AS09, Beyn12}, i.e.,
\begin{equation}\label{eq-L}
L \ge \max_{k=1, \cdots, n_\mathcal {C}} \left ( \sum^{\eta_k}_{l=1} \mu_{l,k}  \right ).
\end{equation}
In practice $U$ can be chosen as a random matrix, and also, the pre-multiplied matrix $U^{H}$ can be replaced by another $L \times n$ matrix of full row rank.

The poles of $F(z)$ can be retrieved by first constructing a rational interpolation of $F$ and then extracting the poles of the approximation \cite{austin2014numerical,austin2015computing}.
Specifically, let $z_0, \dots, z_{N-1}$ be a set of interpolating points within $\mathcal {C}$ that are not the poles of $F(z)$. Consider a rational function of type $(\mu, \nu)$ that satisfies
\begin{equation}\label{eq-rational-qp}
{P(z_i) \over q(z_i)} = F(z_i), \quad i=0, \dots, N-1,
\end{equation}
where $P(z)$ is a matrix polynomial of degree $\mu$, $q$ is a polynomial of degree $\nu$, and $\mu + \nu = N-1$. Obviously, $\nu$ should be larger than or equal to the number of mutually different poles. If $\nu$ is taken to be the number of poles of $F(z)$, Saff \cite{saff1972extension} showed that: (1) for sufficiently large $\mu$ there exist rational interpolants of type $(\mu, \nu)$ that converge uniformly to $F(z)$ in some region; and (2) the poles of the rational interpolants converge to those of $F(z)$. These results offer a possibility to extract the poles of $F$ using the RIA.

In this paper, we follow Jacobi's method that works with the linearized form of interpolation condition \eqref{eq-rational-qp},
\begin{equation}\label{eq-rational-pqf}
P(z_i)  = F(z_i)q(z_i), \quad i=0, \dots, N-1.
\end{equation}
Jacobi's method extracts the poles of a rational function from its moments; see \cite{eugeciouglu1989fast} for the details. In the following, we first extend the idea to the matrix-valued function $F(z)$ and show the relations between the poles and moments of $F(z)$. Then, we show how to compute the eigenvalues and eigenvectors of the original NEP \eqref{eq-gnep} from the moments of $F(z)$.

\subsubsection{Moments of $F(z)$} \label{S-S-momf}  In light of Jacobi's method, the following
moments of $F(z)$ are defined,
\begin{equation}\label{eq-Aalpha-rat}
A_{\alpha} = \sum^{N-1}_{i=0} \omega_i z_i^\alpha F(z_i), \quad \alpha = 0, \dots, 2K-1,
\end{equation}
where $K$ is a finite positive integer such that $2K<N$, $\omega_i$ are the \emph{barycentric weights} for the polynomial interpolations using points $z_i$,
\begin{equation}\label{eq-bary-wt}
\omega_i = {1 \over \prod_{j=0,j \not= i}^{N-1} (z_i - z_j)}.
\end{equation}
Note that in Jacobi's method the highest order of moments is determined by the degree of the denominator polynomial $q(z)$, i.e., $\nu$. However, here $K$ is used in lieu of $\nu$ in order to be consistent with the notations of the CIA in
\cite{AS09, Beyn12}.

To see the connection between the moments $A_\alpha$ and the eigenvalues of $T(z)$, the relations \eqref{eq-thm1-mf} and \eqref{eq-F} are inserted into the moment expression \eqref{eq-Aalpha-rat},
\begin{equation}\label{eq-Aalpha-ratall}
A_{\alpha}  = \sum^{N-1}_{i=0}  \omega_i z_i^\alpha   \left[ U^H V_\mathcal{C}  \Phi(z_i)  W_\mathcal{C}^H   U \right] + \sum^{N-1}_{i=0} \omega_i z_i^\alpha \left[U^H R_\mathcal{C}(z_i) U\right].
\end{equation}
The summation about $R_\mathcal{C}$ accounts for the influence of the eigenvalues outside the contour $\mathcal {C}$. To be able to accurately extract all the eigenvalues within $\mathcal {C}$, this term has to be effectively removed. We attain this goal by invoking the polynomial interpolation theory. In fact, due to the definition of $\omega_i$ \eqref{eq-bary-wt},
the summation gives the leading coefficient matrix of the matrix polynomial of degree at most $N-1$ that interpolates the matrix-valued function $z^\alpha \left[U^H R_\mathcal{C}(z) U\right]$ at the node points $z_i, \,i=0, \dots, N-1$; see \cite{eugeciouglu1989fast}, thus it becomes zero if $z^\alpha \left[U^H R_\mathcal{C}(z) U\right]$ is a polynomial of degree strictly less than $N-1-\alpha$. In practice, since $R_\mathcal{C}(z)$ is analytic, it can be well approximated by polynomials of degrees high enough, which means that it holds
\begin{equation}\label{eq-filterR}
\sum^{N-1}_{i=0} \omega_i z_i^\alpha \left[U^H R_\mathcal{C}(z_i)U\right]  \rightarrow 0, \quad \alpha = 0, \dots, 2K-1\, \mathrm{and}\, N\rightarrow \infty.
\end{equation}
The magnitude of the summation depends on the distribution of $z_i$, $N$ and $\alpha$. In the following, we assume that $N$ and $K$ are chosen such that the magnitude of \eqref{eq-filterR}
is negligible. The criterions for the selection of $N$ and $K$ will be stated in Section \ref{S-S-rialgorithm}. The above assumption leads to the following approximation to the moments
\begin{equation}\label{eq-Aalpha-rat0}
A_{\alpha}  \approx   U^H V_{\mathcal {C}}  \left[\sum^{N-1}_{i=0}  \omega_i z_i^\alpha\Phi(z_i)\right]  W_{\mathcal {C}}^H   U.
\end{equation}

Now we consider the summation about the matrices $\Phi(z_i)$ in \eqref{eq-Aalpha-rat0}. The summation about the matrices can be performed element-wisely, thus it is sufficient to consider the following weighted summation of a nonzero element of $\Phi(z)$ \eqref{eq-MPhiz},
\begin{equation}\label{eq-faj}
\phi_{\alpha,j}(z) := \sum^{N-1}_{i=0} {\omega_i z_i^\alpha \over (z_i - z)^{j} }, \quad j=1,2,\cdots.
\end{equation}
The above definition implies that $\phi_{\alpha,j}(z), \, j>1$ can by expressed by $\phi_{\alpha,1}(z)$,
\begin{equation}\label{eq-gs1}
\phi_{\alpha,j}(z) = {1 \over (j-1)! } \phi_{\alpha,1}^{(j-1)}(z).
\end{equation}
Since we require that $\alpha < 2K<N$, the monomial $z^\alpha$ can be exactly recovered from its values at the $N$ node points $z_i$ by using the barycentric formula of polynomial interpolation theory \cite{berrut2004barycentric}, leading to a close form expression for $\phi_{\alpha,1}(z)$,
\begin{equation}\label{eq-fa1}
\phi_{\alpha,1}(z) = \sum^{N-1}_{i=0} {\omega_i z_i^\alpha \over z_i - z }  = - { z^\alpha \over l(z)},
\end{equation}
where $l(z) = (z-z_0)(z-z_1)\cdots(z-z_{N-1})$ is the node polynomial. It follows from \eqref{eq-fa1} that
\begin{equation}\label{eq-gs2}
z \phi_{\alpha,1}(z) = \phi_{\alpha+1,1}(z).
\end{equation}
By successively differentiating both sides of \eqref{eq-gs2} with respect to $z$ and invoking \eqref{eq-gs1}, one obtains
\begin{equation}\label{eq-gs3}
z \phi_{\alpha,j}(z) + \phi_{\alpha,j-1}(z) = \phi_{\alpha+1,j}(z).
\end{equation}

Let $\Phi_{\alpha,d}$ and $J_{\lambda,d}$ be the $d$-dimensional upper triangular matrix and the $d$-dimensional Jordan block with diagonal entries $\lambda$ defined as follows
$$
\Phi_{\alpha,d}(\lambda) = \begin{bmatrix}
           \phi_{\alpha,1}(\lambda)  & \phi_{\alpha,2}(\lambda) & \cdots  &  \phi_{\alpha,d}(\lambda) \\
             & \ddots & \ddots &  \vdots \\
             &   & \phi_{\alpha,1}(\lambda) & \phi_{\alpha,2}(\lambda) \\
             &   &   & \phi_{\alpha,1}(\lambda)  \\
         \end{bmatrix},
\quad
J_{\lambda, d} = \begin{bmatrix}
           \lambda  & 1 &   &   \\
             & \ddots & \ddots &   \\
             &   & \lambda & 1 \\
             &   &   & \lambda  \\
         \end{bmatrix} \in \mathbb{C}^{d \times d}.
$$
Then, it can be deduced from \eqref{eq-gs2} and \eqref{eq-gs3} that
\begin{equation}\label{eq-Ja}
\Phi_{\alpha,d}(\lambda) = J_{\lambda,d} \Phi_{\alpha-1,d}(\lambda) = \cdots  = J_{\lambda,d}^\alpha \Phi_{0,d}(\lambda).
\end{equation}

Applying \eqref{eq-Ja} to each diagonal block of $\Phi(z_i)$ in the summation \eqref{eq-Aalpha-rat0}, one obtains
\begin{equation}\label{eq-Aalpha-brk}
\sum^{N-1}_{i=0}  \omega_i z_i^\alpha\Phi(z_i) = \Lambda^\alpha \Phi,
\end{equation}
where the matrix $\Lambda$ has Jordan normal form
\begin{equation} \label{eq-Mpr-Jordan}
\Lambda = \begin{bmatrix}
  J_1 &   &   \\
    & \ddots &   \\
    &   & J_{n_\mathcal {C}} \\
\end{bmatrix}, \quad
J_k = \begin{bmatrix}
  J_{k}^{1} &   &   \\
    & \ddots &   \\
    &   & J_{k}^{\eta_k} \\
\end{bmatrix}, \quad
J_{k}^{l} = J_{\lambda_k, \mu_{l,k} } \in \mathbb{C}^{\mu_{l,k} \times \mu_{l,k}},
\end{equation}
and the constant matrix $\Phi$ has the same structure as $\Phi(z)$ in \eqref{eq-MPhiz},
\begin{equation} \label{eq-Mpr-Jordan-ci}
\Phi = \begin{bmatrix}
  \Phi_1 &   &   \\
    & \ddots &   \\
    &   & \Phi_{n_\mathcal {C}} \\
\end{bmatrix}, \quad
\Phi_k = \begin{bmatrix}
  \Phi_{k}^{1} &   &   \\
    & \ddots &   \\
    &   & \Phi_{k}^{\eta_k} \\
\end{bmatrix}, \quad
\Phi_{k}^{l} = \Phi_{0, \mu_{l,k} }(\lambda_k).
\end{equation}
Note that $\Phi$ is an upper triangular matrix whose diagonal entries are given by $\phi_{0,1}(\lambda_k) = -  1 / l(\lambda_k), \, k = 1, \dots, n_\mathcal {C}$. Since we assume that the
sampling points $z_i$ do not coincide with the poles $\lambda_k$, $l(\lambda_k)\not = 0$ and thus $\Phi$ is nonsingular.

Inserting \eqref{eq-Aalpha-brk} into \eqref{eq-Ja} leads to the following explicit relation between the moments $A_\alpha$ and the eigenvalues of $T(z)$,
\begin{equation}\label{eq-Aalpha-rat-mat}
A_{\alpha}  \approx U^H V_\mathcal {C} \Lambda^\alpha \Phi  W_\mathcal {C}^H U.
\end{equation}

\subsubsection{Computing eigen-solutions} \label{S-S-eigsol}

To extract $\Lambda$ from $A_\alpha$, we define two block Hankel matrices,
\begin{equation} \label{eq-HH}
H:= \begin{bmatrix}
      A_0     & A_1 & \cdots & A_{K-1} \\
      A_1     & \ddots     & \ddots &         \vdots \\
      \vdots         & \ddots     & \ddots &          \vdots \\
      A_{K-1} & \cdots     & \cdots & A_{2K-2} \\
    \end{bmatrix}, \quad
H^{<}:= \begin{bmatrix}
      A_1       & A_2 & \cdots & A_{K} \\
      A_2       & \ddots     & \ddots &           \vdots \\
      \vdots           & \ddots     & \ddots &           \vdots \\
      A_{K} & \cdots     & \cdots & A_{2K-1} \\
    \end{bmatrix}.
 \end{equation}
Then the poles of $F(z)$ within $\mathcal {C}$ can be approximated by the eigenvalues of the generalized eigenvalue problem for the pencil $H^< - \lambda H$. In fact, it follows from \eqref{eq-Aalpha-rat-mat} that
\begin{equation}\label{eq-HHVW}
H \approx V_{[K]} W_{[K]}^H  \quad \mathrm{and} \quad
H^{<} \approx V_{[K]} \Lambda W_{[K]}^H,
\end{equation}
where
\begin{equation} \label{eq-VWK}
V_{[K]}= \begin{bmatrix}
      U^H V_\mathcal {C} \\
      U^H V_\mathcal {C}\Lambda \\
      \vdots \\
      U^H V_\mathcal {C}\Lambda^{K-1} \\
    \end{bmatrix}
 \quad \mathrm{and} \quad
W_{[K]}^H =
    \begin{bmatrix}
      \Phi  W_\mathcal {C}^H U & \Lambda  \Phi W_\mathcal {C}^H U & \cdots & \Lambda^{K-1}  \Phi W_\mathcal {C}^H U \\
    \end{bmatrix}.
\end{equation}
Obviously, the poles $\lambda_k,\, k=1, \dots \bar n_{\mathcal {C}}$ coincide with the eigenvalues of the pencil $H^< - \lambda H$ if $K$ is chosen to fulfill the following rank condition,
\begin{equation} \label{eq-ranks}
\rank (V_{[K]})
    =
\rank (W_{[K]}^H)
    = \bar n_{\mathcal {C}},
\end{equation}
which requires
$K\cdot L \ge \bar n_{\mathcal {C}}$.

Till now, the problem of finding the eigenvalues of $T(z)$ has been transformed to the solution of the generalized eigenvalue problem about the pencil $H^< - \lambda H$.
Following the procedures in \cite{Beyn12}, we first compute the truncated SVD of $H$,
\begin{equation} \label{eq-svdH}
 H \approx V_0 \Sigma_0 W_0^H,
\end{equation}
where $V_0$ and $W_0$ are matrices with orthogonal columns, and the number of singular values retained in $\Sigma_0$ determines the number of eigenvalues $\bar n_{\mathcal {C}}$, which will be
detailed in Section \ref{S-S-rialgorithm}. Then we compute the matrix
\begin{equation} \label{eq-A}
A = V_0^H H^{<} W_0 \Sigma_0^{-1},
\end{equation}
which has Jordan normal form $\Lambda$. The eigenvalue problem for $A$ can be easily solved using the Matlab function `\verb"eig"'. Given the eigenvalue decomposition $A=B \Lambda B^{-1}$ with $B$ being the eigenvector matrix, the matrices $V_{[K]}$ and $W_{[K]}^H$ in \eqref{eq-HHVW} can then be computed as $V_{[K]} = V_0 B$ and $W_{[K]}^H = B^{-1} \Sigma_0 W_0^H$. To compute the eigenvectors of the original NEP \eqref{eq-gnep}, one needs to formulate a new matrix $M$ as follows
\begin{equation}\label{eq-Mspace}
M = [M_0\; M_1\; \cdots \; M_{K-1} ],
\end{equation}
where
\begin{equation}\label{eq-Malpha}
\begin{split}
M_{\alpha} = \sum_{i=0}^{N-1} \omega_i z_i^\alpha T(z_i)^{-1} U, \quad \alpha=0,1,\dots, K-1.
\end{split}
\end{equation}
Note that the definition of $M_{\alpha}$ is similar to that of $A_\alpha$ in \eqref{eq-Aalpha-rat}, thus $M_{\alpha}$ can be obtained as by-products during the computation of $A_\alpha$.
The matrix $V_\mathcal {C}$ relates to $M$ via
$$
M \approx V_\mathcal {C} W_{[K]}^H = V_\mathcal {C} B^{-1} \Sigma_0 W_0^H \quad \Rightarrow \quad V_\mathcal {C} \approx M W_0 \Sigma_0^{-1} B.
$$
Hence,  if $(b, \lambda')$ is an eigenpair of $A$, the corresponding eigenpair $(v, \lambda)$ of the original NEP \eqref{eq-gnep} can be obtained as
\begin{equation} \label{eq-eigens}
\lambda \approx \lambda' \quad \mathrm{and} \quad v \approx M W_0 \Sigma_0^{-1} b.
\end{equation}

So far, we have described the RIA for solving general NEPs. By employing Jacob's method \cite{eugeciouglu1989fast} and introducing the moments $A_\alpha$ of the resolvent $T(z)^{-1}$ \eqref{eq-Aalpha-rat}, we obtain a generalized eigenvalue problem for pencil $H^< - \lambda H$, whose eigenvalues approximate the eigenvalues of the original NEP. The accuracy of the eigenvalue and eigenvector approximations is depend on how well the analytic function $R_\mathcal{C}(z)$ is approximated by the polynomial interpolation associated with the points $z_i$ and barycentric weights $\omega_i$; see \eqref{eq-filterR}. Generally, for a given value of $K$ that satisfying the rank condition \eqref{eq-ranks}, increasing the number $N$ of sampling points $z_i$ will improve the quality of eigenpair approximations \eqref{eq-eigens}.

\subsection{Numerical algorithm and practical considerations} \label{S-S-rialgorithm}

The RIA in Section \ref{S-S-ria} directly leads to the following procedures for solving NEPs:

\begin{enumerate}
  \item Initialization. Fix the contour $\mathcal {C}$, the number $N$ and the sampling points $z_i, i=0,1,\dots, N-1$ on or within the contour; compute the corresponding weights $\omega_i$; fix the number $L$ and generate a $n \times L$ random matrix $U$; choose a suitable value of $K$. \label{alg-ss-step1}

  \item Compute $T(z_i)^{-1} U$ for $i=0,1,\dots, N-1$. \label{alg-ss-step2}

  \item Compute moments $A_\alpha,\, \alpha = 0, \dots, 2K-1$ from \eqref{eq-Aalpha-rat}, and $M_\alpha, \, \alpha=0,1,\dots, K-1$ from \eqref{eq-Malpha}. \label{alg-ss-step3}

  \item Construct two block Hankel matrices $H$ and $H^<$ according to \eqref{eq-HH}; construct $M$ according to \eqref{eq-Mspace}. \label{alg-ss-step4}%

  \item Compute the truncated SVD $H\approx V_0 \Sigma_0 W_0^H$ and determine the number of eigenvalues $\bar n_{\mathcal {C}}$. \label{alg-ss-step5}

  \item Form the matrix $A$ according to \eqref{eq-A}, and solve the standard eigenvalue problem for $A$. Let $(b, \lambda)$ denote an eigenpair of $A$, then the corresponding eigenpair of the original NEP \eqref{eq-gnep} is given by $(v, M W_0 \Sigma_0^{-1} b)$. \label{alg-ss-step6}
\end{enumerate}

This algorithm will be denoted by \emph{SS-RI}. In the RIA, the sampling points $z_i$ can be a general set of interpolation points within or on the contour $\mathcal {C}$. When $z_i$ are chosen as the nodes of some numerical quadrature rules on the contour $\mathcal{C}$, SS-RI becomes SS-CI in \cite{AS09}. In addition, it has been proved in \cite{austin2014numerical, austin2015computing} that for linear eigenvalue problems the CIA and RIA are indeed mathematically equivalent (in exact arithmetic) when $z_i$ uniformly distributed on a circle.

In the following some practical issues of the SS-RI algorithm are discussed.

\begin{enumerate}
  \item{\it How to choose sampling points $z_i$.} The sampling points serve as the interpolation points. They should be chosen so that the residual matrix $R_\mathcal{C}(z)$ can be well approximated by a polynomial of the lowest degree. In practice, the domain of interest is often an interval or a rectangle, then the sampling points can be set as the Chebyshev points in the domain or the quadrature points on the boundary contour.

  \item{\it How to determine $L$, $K$ and $N$.} Theoretically, $L$ only has to be larger than or equal to the maximal algebraic multiplicity of eigenvalues of $T(z)$ in $\mathcal {C}$. However, in practical applications of the SS-RI algorithms $L$ would have to be more larger \cite{sakurai2013efficient}. $K$ should be chosen such that size of the search space $M$ is large enough, i.e., $K\cdot L \ge \bar n_{\mathcal {C}}$. On the contrary, for a fixed $N$ a larger $K$ may cause accuracy decline because it diminishes the performance of the interpolation \eqref{eq-filterR}. In practice, $K \leqslant N/4$ can be used \cite{sakurai2013efficient}. The choice of $N$ is mainly determined by the tradeoff between the accuracy requirement and the computational cost.

  \item{\it How to avoid under- and overflow in computing weights using \eqref{eq-bary-wt}.}
A strategy is to simplify the weights by cancelling the common factors. This operation clearly does not alter the eigenvalues of the pencil $H^< - \lambda H$. For certain special sets of interpolation points, one can give explicit formulas for the simplified barycentric weights \cite{berrut2004barycentric, wang2014explicit}. A typical example is the Chebyshev points of the first kind on the unit interval $[-1,1]$,
\begin{equation} \label{eq-Chebyshev-pts}
x_i = \cos { (2i +1) \pi\over 2N}, \quad i = 0, \dots, N-1.
\end{equation}
In this case after cancelling factors independent of $i$ one finds
\begin{equation} \label{eq-Chebyshev-wts}
\omega_i = (-1)^i \sin { (2i +1) \pi\over 2N},
\end{equation}
which vary by factors $O(N)$. For more general situations, one can multiply each factor $z_i - z_j$ in \eqref{eq-bary-wt} by a scaling coefficient $C^{-1}$, with $C$ being the capacity of the interval concerned; see \cite{berrut2004barycentric} for the details.

  \item{\it How to determine the total number of eigenvalues $\bar n_{\mathcal {C}}$.} From \eqref{eq-svdH} one knows that $\bar n_{\mathcal {C}}$ is given by the number of leading singular values retained in $\Sigma_0$. Hence, one can determine it by detecting the largest gap
in the singular values of $H$ in \eqref{eq-svdH}. In case $N$, $L$ and $K$ are properly chosen and $z_i$ do not coincide with the eigenvalues, the largest gap between two successive singular values denoted by
$$
g_{\max} = \max_{j=1,\cdots, K\cdot L - 1}\left({\sigma_j / \sigma_{j+1}}\right),
$$
should reach its maximum at $j=\bar n_{\mathcal {C}}$. This leads to the following strategy,
\begin{equation} \label{eq-gap}
\mathrm{if} \;  g_{\max} \ge \tol_{\gap},  \quad \mathrm{then}\; \bar n_{\mathcal {C}} = \underset{ j=1,\cdots, K\cdot L - 1 }{ \mathrm{Argmax} }  \left({\sigma_j / \sigma_{j+1}} \right),
\end{equation}
where the user-defined constant $\tol_{\gap}$ is used only to check whether a $g_{\max}$ is reasonable; in this work, we always set $\tol_{\gap} = 10^{3}$. If condition $ g_{\max} \ge \tol_{\gap}$ is not satisfied, the obtained $\bar n_{\mathcal {C}}$ may be not correct, and thus the computed eigenpairs may be wrong or of low accuracy. In this case, checking the residuals of the eigenpairs, $||T(\lambda)v||_2/||v||_2$, would be helpful. Usually, improvement to the results can be achieved by increasing $N$, $L$ or $K$.

Similar strategies have been used in \cite{AS09, Beyn12}. But in those papers $\bar n_{\mathcal {C}}$ is determined by the truncation of the singular values of $H$
using a predetermined threshold $tol$. In general the optimal threshold $tol$ is affected by several factors, including the conditioning of the matrices in \eqref{eq-ranks}, the distribution of the eigenvalues and the positions of the contour, and thus may not be easily determined before solving the NEPs.

\end{enumerate}

Let us close this section with a summary of the proposed RIA and the SS-RI algorithm. The RIA is a generalization of the CIA in the sense that polynomial interpolation is used to extract the poles of rational functions instead of contour quadrature. This generalization explains the phenomenon described in Section \ref{S-S-motivation}, and offers new possibilities to improve the accuracy and reduce the computational cost.
The SS-RI algorithm is an straightforward implementation of the RIA. In Section \ref{S-S-nin-ss}, we will show by a numerical example that in general SS-RI and SS-CI can achieve the same accuracy, but for real NEPs SS-RI can save a half of the CPU time and memory by using real sampling points and performing real arithmetic.
The fact that the RIA allows us to improve the accuracy of eigen-solutions will be demonstrated by using the \ssrrs{} algorithm in the following section.

The SS-RI algorithm is robust and accurate if a large $L$ but a small $K$ are used. However, for large-scale problems, a small $L$ is essential to reduce the computational burden. To decrease $L$, one has to increase $K$ due to the rank condition $K\cdot L \ge \bar n_{\mathcal {C}}$ \eqref{eq-ranks}, but this would make the two matrices $H$ and $H^{<}$ rank-deficient and finally make the algorithm unstable and inaccurate; see Section \ref{S-S-nin-ss} for the related numerical results.
Therefore, the SS-RI algorithm is suitable for solving small-scale NEPs. For this class of problems, the algorithm can even be made more concise by discarding the matrix $U$ in \eqref{eq-F} and directly using $F(z) = T^{-1}(z)$. This corresponds to using the largest $L$, thus the algorithm could achieve its optimal robustness. In this case, usually only a small value of $K$ (e.g., $K<10$) is needed to compensate for the possible rank deficiency of $V_\mathcal {C}$ and $W_\mathcal {C}$. This concise version of the SS-RI algorithm will be referred to as the \emph{SS-FULL} and will be used to solve the projected NEP by the Rayleigh-Ritz procedure.

\section{Rayleigh-Ritz reformulation} \label{S-Rayleigh}

As mentioned before, the SS-RI algorithm would become unstable and inaccurate when a small value of $L$ is used.
This problem with can be effectively overcome by using the Rayleigh-Ritz projection procedure \cite{YS13}.
In this section, we describe the sampling scheme to construct approximate eigenspaces based on the RIA and the RSRR algorithm for solving general NEPs. The RSRR algorithm has good stability and accuracy when a reasonably small value of $L$ is used, and is suitable for large-scale problems.

The classical Rayleigh-Ritz procedure relies on a proper search space that contains the interested eigenvectors. Once such a search space is available and let $Q\in \mathbb{C}^{n \times k}$ be an orthogonal basis of it, then the original NEP \eqref{eq-gnep} can be converted to the following reduced NEP
\begin{equation} \label{eq-Ts-NEP}
T_Q(\lambda)g = 0 \quad \mathrm{with} \quad T_Q(z) = Q^H T(z)Q \in \mathbb{C}^{k \times k}.
\end{equation}
Let $(\lambda, g)$ be any eigenpair of the reduced NEP, then $(\lambda, Sg)$ is an eigenpair of the original NEP.

For ease of presentation, in the following we first describe the moment scheme based on the new RIA, and then introduce the sampling scheme from the moment scheme.

\subsection{Resolvent moment scheme} \label{S-S-moment-ria}

Let $z_i, i=0, \cdots, N-1$ be a set of sampling (interpolation) points in the domain enclosed by the contour $\mathcal {C}$, and denote by $\omega_i$ the corresponding barycentric weights \eqref{eq-bary-wt}.
We begin by defining the moments $M_\alpha$ of $T^{-1}(z)U$ up to the order $K-1$ ($K\leqslant N$) as \eqref{eq-Malpha}. Then, a similar derivation to the relation \eqref{eq-Aalpha-rat-mat} leads to
\begin{equation}\label{eq-mom-apr}
M_{\alpha}  \approx  V_\mathcal {C} \Lambda^\alpha \Phi  W_\mathcal {C}^H U, \quad \alpha = 0, \dots, K-1\, \mathrm{and}\, N\rightarrow\infty.
\end{equation}
The key to the approximation \eqref{eq-mom-apr} is the effective elimination of $R_\mathcal{C}(z)$ in the expression of the resolvent $T(z)^{-1}$ \eqref{eq-thm1-mf}. In fact, due to the fact that $R_\mathcal{C}(z)$ is analytic and $\omega_i$ are the barycentric weights associated with $z_i$, a similar argument to \eqref{eq-filterR} leads to
\begin{equation}\label{eq-rmom}
\sum^{N-1}_{i=0} \omega_i z_i^\alpha R_\mathcal{C}(z_i)U  \rightarrow 0, \quad \alpha = 0, \dots, K-1\, \mathrm{and}\, N\rightarrow\infty.
\end{equation}
Note that here only the first $K$ moments of $R_\mathcal{C}(z)U$ are required to be negligible, but in \eqref{eq-filterR} the order is up to $2K$.

To retrieve the eigenspace $\spanset (V_{\mathcal{C}})$, we collect all the $K$ moments together,
\begin{equation}\label{eq-M-pr}
\begin{split}
M &:= \left ( M_0, \, M_1, \, \cdots, \,  M_{K-1} \right) \\
& \approx  V_{\mathcal{C}}
\begin{bmatrix}
      \Phi  W^H_{\mathcal{C}} U & \Lambda \Phi  W^H_{\mathcal{C}} U & \cdots & \Lambda^{K-1} \Phi  W^H_{\mathcal{C}} U\\
    \end{bmatrix} \\
    &= V_{\mathcal{C}} W_{[K]}^H.
\end{split}
\end{equation}
It follows from expression \eqref{eq-M-pr} that if $K$ is chosen such that
\begin{equation}\label{eq-rankcond-M}
K\cdot L  \geqslant \rank \left(W_{[K]}^H \right)  \geqslant  \rank (V_{\mathcal{C}}),
\end{equation}
then the matrix $M$ will span approximately the same subspace as $V_{\mathcal{C}}$, which means that
there exist a positive integer $K_0$ such that
\begin{equation}\label{eq-space-MV}
K_0< K\leqslant N\rightarrow \infty \quad \Rightarrow \quad \spanset (M)  = \spanset(V_{\mathcal{C}}).
\end{equation}
For a finite number $N$ the discrepancy of the two subspaces diminishes with the increase of $N$. Similar conclusion has been obtained for the CIA; see Theorem 3 in \cite{YS13}.

It turns out that the eigenspace of the moment scheme is the same as the one used in the SS-RI algorithm; see \eqref{eq-eigens}, \eqref{eq-Mspace} and \eqref{eq-Malpha}. However the derivation here is directly from Theorem \ref{thm1} and the interpolation theory; one does not need to refer to the derivation of the SS-RI algorithm in Section \ref{S-S-ria}. In addition, the above procedure for constructing eigenspaces is similar to those based on the CIA \cite{YS13}, but the RIA offers more freedom to the distribution of the sampling points.

 Finally, we remark that in practice the eigenspaces generated by using the moment scheme would be unreliable for two reasons. First, the moment matrix $M$ \eqref{eq-M-pr} tends to be rank-deficient for large $K$, and finally the rank condition \eqref{eq-rankcond-M} would not be fulfilled securely (see Section \ref{S-S-nin-ssrim} for numerical evidences), and second, improper computational procedures for $\omega_i$ would often cause under- and overflow (see Section \ref{S-S-rialgorithm}). These problems will be completely avoided by the sampling scheme in the following section.

\subsection{Resolvent sampling scheme} \label{S-S-samp}

The sampling scheme is inspired by the fact that each $M_\alpha$ \eqref{eq-Malpha} can be seen as a combination of $T(z_i)^{-1} U$, i.e.,
\begin{equation*}\label{eq-malpha-lin}
M_\alpha = \sum^{N-1}_{i=0} c_{\alpha, i} T^{-1}(z_i)U, \quad \alpha=0,1,\dots, K-1,
\end{equation*}
where, $c_{\alpha, i} = \omega_i z_i^\alpha$ are the combining coefficients determined by $\alpha$ and $z_i$. Although the above combination is indeed nonlinear, for a given set of points $z_i$ the columns of $T^{-1}(z_i)U$ should span a larger subspace than $M_\alpha$. Therefore, we propose to generate eigenspaces directly using the former.
Following this idea, we form the following sampling matrix by collecting all $T(z_i)^{-1} U$,
\begin{equation} \label{eq-hatS}
S = \left [ T(z_0)^{-1} U, \, T(z_1)^{-1} U, \, \cdots, \, T(z_{N-1})^{-1} U \right] \in \mathbb{C}^{n \times N\cdot L},
\end{equation}
and use $\spanset (S)$ as the eigenspace in the sampling scheme.
The above argument indicates that the proposed sampling scheme actually constructs a larger subspace, i.e.,
\begin{equation} \label{eq-ms}
\spanset (M) \subseteq \spanset (S).
\end{equation}
The equality holds when $K=N$ and $M$ is not rank-deficient.

When used in the Rayleigh-Ritz projection, it is important for the constructed subspaces to contain all the eigenvectors. The sampling scheme can meet this requirement.
Theoretically, when $N$ and $K$ are chosen such that the relation \eqref{eq-space-MV} holds, then it follows from \eqref{eq-ms} that
\begin{equation} \label{eq-SV}
\spanset (V_{\mathcal{C}}) \subseteq  \spanset (S), \quad N\rightarrow \infty,
\end{equation}
In practice, $N$ is a finite number and thus the influence of $R_\mathcal{C}(z)$ \eqref{eq-rmom} may not be completely removed. However, according to \eqref{eq-SV} the
subspaces of $S$ will approximately contain the eigenspaces when $N$ is large enough. Note that the following condition is a prerequisite
\begin{equation} \label{eq-rank-NL}
      N\cdot L \ge  \rank(S) \ge \rank(V_{\mathcal{C}}).
\end{equation}

It remains to show that the sampling scheme is more robust than the moment scheme. This indeed can be seen from the rank property of the matrices $S$ and $M$. According to the definition of the moments \eqref{eq-Malpha}, one has
\begin{equation} \label{eq-MS}
M = S Z,
\end{equation}
where $Z$ is a Vandermonde-like matrix
$$
Z = \begin{bmatrix}
           Z_0^{[0]}       & Z_0^{[1]}         & \cdots  &  Z_0^{[K-1]} \\
           Z_1^{[0]}       & Z_1^{[1]}         & \cdots  &  Z_1^{[K-1]} \\
                     \vdots      & \vdots  & \ddots  & \vdots \\
           Z_{N-1}^{[0]}       & Z_{N-1}^{[1]}     & \cdots  &  Z_{N-1}^{[K-1]}  \\
         \end{bmatrix} \in \mathbb{C}^{N\cdot L \times K\cdot L}
         \quad \mathrm{and} \quad Z_i^{[\alpha]} = \begin{bmatrix}
           \omega_i z_i^\alpha       &          &   \\
                  & \ddots         &   \\
                 &       &  \omega_i z_i^\alpha  \\
         \end{bmatrix}
        \in \mathbb{C}^{L \times L}.
$$
Obviously, $Z$ would become rank-deficient when $K$ is large, which accounts for the possible rank-deficiency of $M$. On the other hand, $S$ should have better rank property than $M$.

In the Rayleigh-Ritz procedure, an orthogonal basis of $S$, denoted by $Q$, is needed. It can be computed by the truncated singular value decomposition (SVD) of $S$ with a tolerance $\delta$.
The computational cost of SVD scales linearly with $n$ but quadratically with the product $N\cdot L$, therefore the reduction of $N\cdot L$ becomes a major concern when the work of SVD becomes dominated. Let $k$ denote the numerical rank of $S$.
In order to properly extract all the eigenvalues inside $\mathcal {C}$, the condition \eqref{eq-rank-NL} has to be satisfied. However, since $\bar n_{\mathcal {C}} \ge \rank(V_{\mathcal{C}})$, it is more convenient to use the following condition in practice
\begin{equation} \label{eq-rank-nc}
N\cdot L \ge k \ge \bar n_{\mathcal {C}},
\end{equation}
which means that $N$ and $L$ have to been chosen such that $k$ is larger than the total number of eigenvalues.
Of course, $L$ should be not less that the maximal algebraic multiplicity of the eigenvalues in $\mathcal {C}$.

There is a situation that deserves special attention in practice. When some sampling points approach to the unknown eigenvalues, the corresponding matrices $T(z_i)$ will become nearly singular, thus $T(z_i)^{-1}U$ will be of very large values and dominated in the matrix $S$. This would enhance the accuracy of the eigenvalues near the sampling points, but on the contrary would considerably deteriorate the accuracy of the other eigenvalues. The way to avoid this issue is to normalize the column vectors of $S$ before the SVD so that they are of the same norm. This operation does not alter the eigenspaces but guarantees the robustness of the sampling scheme.

\subsection{The \ssrrs{} algorithm} \label{S-S-algorithm}

The main eigensolver \ssrrs{}, the resolvent sampling based Rayleigh-Ritz method, consists of the following procedures:
\begin{enumerate}
  \item Initialization. Fix the contour $\mathcal {C}$, the number $N$ and the sampling points $z_i, i=0,1,\dots, N-1$ on or within the contour. Fix the number $L$ and generate a $n \times L$ random matrix $U$. \label{alg-cirrip-step1}

  \item Compute $T(z_i)^{-1} U$ for $i=0,1,\dots, N-1$. \label{alg-cirrip-step2}

  \item Form $S$ according to \eqref{eq-hatS}. Generate the matrix $Q$ via the truncated SVD $S \approx Q \Sigma V^H $, where the first $k$ singular values larger than $\delta\cdot\sigma_1$ is retained. \label{alg-cirrip-tsvd}

  \item Compute $T_Q(z) = Q^H T(z)Q$, and solve the projected NEP $T_Q(\lambda)g = 0$ using the SS-FULL algorithm in Section \ref{S-S-rialgorithm} to obtain $\bar n_{\mathcal {C}}$ eigenpairs $(g_j, \lambda_j), \, j=1, \cdots, \bar n_{\mathcal {C}}$. \label{alg-cirrip-step4}

  \item Compute the eigenpairs $(v_j,\lambda_j), \, j=1, \cdots, \bar n_{\mathcal {C}}$ of the original NEP \eqref{eq-gnep} via $v_j = Q g_j$; the eigenvalues of the original NEP are equal to those of the projected NEP. \label{alg-cirrip-step5}
\end{enumerate}

In Step \ref{alg-cirrip-step1}, the foremost task is to choose a suitable $\mathcal {C}$ in which the eigenvalues are searched for. For many engineering problems, one may have a-priori information about the locations of the interested eigenvalues from experiments, theoretical predictions, eigenvalues of similar structures, etc.
Another cheaper way to acquire such information is to perform an analysis on a coarse mesh or a simplified model, both of which have been frequently used in literature; see e.g., \cite{Voss07}.
The contour is usually chosen to be an \emph{ellipse} or a \emph{rectangule} enclosing the interested eigenvalues. The sampling points can then be set as the quadrature points on the contour or a special set of points, e.g., the Chebyshev points, within the contour.

Given a value of $L$, the number $N$ should be chosen such that the ratio of the first to the last singular values of $S$ is large enough, e.g., larger than $10^{14}$. This is often the case when $N \cdot L$ is 2 or 3 times larger than the number of the eigenvalues $\bar n_{\mathcal {C}}$; see \eqref{eq-rank-nc}. The number $L$ has to be at least not less than the algebraic multiplicity of the interested eigenvalues. In some cases (e.g., structural modal analysis) this can be estimated from the symmetry of the problems. Otherwise, one can use a slightly large $L$.

Concerning the selection of $N$ and $L$, one should also consider the computational costs and accuracy of the \ssrrs{} algorithm. The effects of $N$ and $L$ on the computational costs lie in two aspects:
the orthogonalization in the sampling scheme as mentioned before, and the solution of the linear system of equations $T(z_i)^{-1} U$.
The computational work related to the latter depend on the linear solvers used.
In the FEM, usually sparse direct solvers are preferable due to its good stability. In this situation, the solution
for problems with multiple right-hand-sides could be cost-effective once a factorization (e.g., the LDU and
Cholesky factorizations) of the system matrix is accomplished. Thus, $N$ determines the main computational cost.
However, in BEM or FEM for too large problems, iterative solvers have to be used.
Then one would pay much more for using a large value of $L$, even with some advanced solvers for multiple right-hand-side problems, e.g., the subspace recycling strategy \cite{parks2006recycling}, because the $L$ right-hand-sides are linearly independent. In this situation,
$N$ and $L$ may have the approximately same weights in determining the computational costs. The relative importance of $N$ and $L$ on the accuracy of \ssrrs{} is so far not clear, and will be numerical studied in Section \ref{S-S-nin-NL}.

In Steps \ref{alg-cirrip-step2} and \ref{alg-cirrip-tsvd} the eigenspaces are constructed using the sampling scheme. The threshold $\delta$ of the truncated SVD can be set as $\delta = 10^{-14}$. In Steps  \ref{alg-cirrip-step4} and \ref{alg-cirrip-step5} the reduced NEP is solved using the contour in Step \ref{alg-cirrip-step1}, and the eigenpairs of the orginal NEP are computed.

The \ssrrs{} algorithm can be easily implemented in conjunction with other software. In fact, only Steps \ref{alg-cirrip-step2} and \ref{alg-cirrip-step4} involve the operations with the matrix $T(z)$ and need to be carried out by the host software. The other steps can be considered as preprocessing and postprocessing parts for the host software. In this paper, \ssrrs{} is implemented into the FEM software \textsc{Ansys}\textregistered{} and our in-house fast BEM code \cite{CWX15}, which will be used to solve the large-scale examples in Section \ref{S-ne-realapps}.

When the eigenspaces are constructed using the moment scheme, one obtains the \ssrrm{} algorithm \cite{YS13}. But, according to the RIA the sampling points are allowed to lie inside the contour.
The \ssrrm{} algorithm can be implemented in a similar way as \ssrrs{}. In this case, the order of moments $K$ and the barycentric weights $\omega_i$ have to be initialized in Step \ref{alg-cirrip-step1}, and in Step \ref{alg-cirrip-tsvd} instead of forming $S$ according to \eqref{eq-hatS} one has to form $M$ from \eqref{eq-M-pr} and \eqref{eq-Malpha}.

Finally, we summarize the two major advantages of the \ssrrs{} algorithm. First, \ssrrs{} is more accurate and robust than \ssrrm{}, because: (1) the sampling scheme generates better subspaces than the moment scheme, and the rank property of the sampling matrix $S$ is better than the moment matrix $M$; and (2) the computations of the moments $M_\alpha$ and barycentric weights $\omega_i$ as required in the moment scheme are completely avoided in the sampling scheme. The latter also leads to a more concise numerical implementation of \ssrrs{}.

The second advantage is that \ssrrs{} allows us to enhance the accuracy of the eigenapirs by optimizing the locations of the sampling points. For example, in Section \ref{S-S-accuracy-gun} we will show that remarkable improvement to the accuracy can be achieved by putting the sampling points close to the eigenvalues. However, in the moment scheme it is commonly known that the sampling points should not approach to any eigenvalue, since this will deteriorate the accuracy of the other eigenpairs \cite{hasegawa2015recovering}.
Note that in this situation the normalization of the column vectors of $S$ is crucial to guarantee the robustness of \ssrrs{}.

\section{Numerical Investigation} \label{S-nin}

In this section, the theoretical benefits of the RIA and the performance of the algorithms SS-RI and \ssrrs{}
are numerically studied and verified.
In Section \ref{S-S-nin-ss} the accuracy and computational costs of SS-RI are studied and compared with SS-CI. In particular, the effects of parameters $L$ and $K$ on the accuracy of SS-RI are tested and the possible loss of accuracy caused by using a large $K$ is illustrated. Sections \ref{S-S-nin-ssrim} to \ref{S-S-NLEIGS-gun} are devoted to the study of the \ssrrs{} algorithm from several aspects:
\begin{itemize}
  \item Section \ref{S-S-nin-ssrim}, compare with SS-RI and \ssrrm{} to show the remarkably improved robustness and accuracy;
  \item Section \ref{S-S-nin-NL}, study the effects of $N$ and $L$ on the accuracy of \ssrrs{};
  \item Section \ref{S-S-accuracy-gun}, confirm that the RIA allows us to improve the accuracy of \ssrrs{} by putting sampling points close to the eigenvalues;
  \item Section \ref{S-S-NLEIGS-gun}, compare with the state-of-the-art nonlinear eigensolver NLEIGS to show the good performance of \ssrrs{}.
\end{itemize}

All the computations were performed on a personal computer with an Intel$^\circledR$
Core$^\mathrm{TM}$ i3-2100 (3.10 GHz) CPU and 16 GB RAM.
The coding and computations are conducted in \textsc{Matlab} R2009a, and the linear systems in computing $T(z_i)^{-1} U$ are solved by using the \textsc{Matlab} backslash operator \verb"`\'" in sparse mode.

\subsection{Numerical investigation of the SS-RI algorithm} \label{S-S-nin-ss}

The performance of SS-RI is investigated by a comparison with SS-CI in \cite{AS09}. We will show that both algorithms can reach the same accuracy, but for real eigenvalue problems the SS-RI algorithm can save almost a half of the CPU time by using real arithmetic. In addition, we will also show the lose of accuracy caused by using a large $K$.

The numerical experiment uses the following example with real eigenvalues.

\begin{example}[\textit{Loaded string}] \label{ex-string}
Consider a special case of the NEP \eqref{eq-REP-string}, where $J=1$, $\sigma_1 = 1$ and $C_j = e_n e_n^T$; that is,
$T(z) =  K_{\rm s} + {z \over z-1} e_n e_n^T - z M$,
where $n$ is the number of equally-spaced finite elements,
$$
K_{\rm s} = {n}\begin{bmatrix}
  2 & -1 &   &   \\
  -1 & \ddots & \ddots &   \\
    & \ddots & 2 & -1 \\
    &   & -1 & 2 \\
\end{bmatrix}
\quad
\mathrm{and}
\quad
M = {1 \over 6n}\begin{bmatrix}
  4 & 1 &   &   \\
  1 & \ddots & \ddots &   \\
    & \ddots & 4 & 1 \\
    &   & 1 & 2 \\
\end{bmatrix}.
$$
The matrix $T(z)$ is symmetric and real, and all the eigenvalues are real \cite{SSI06}.
The $32$ eigenvalues in real interval $\mathcal {I} = [3, 10000]$ are sought. The size of this problem is set to be $n=400$.

The contour for the SS-CI algorithm is an ellipse
defined by
\begin{equation} \label{eq-ellip-contour}
\varphi(\alpha) = \gamma + \left [ a \cos(\alpha) + \im  b \sin(\alpha) \right], \quad \alpha \in [0, 2\pi),
\end{equation}
whose major axis is coincident with $\mathcal {I}$, thus the center is $\gamma = 5001.5$ and the length of the semi-major axis is $a = 4998.5$. The length of the semi-minor axis is set as $b = a/2$.
The contour moments are computed by using the following
$N$-point trapezoidal rule
\begin{equation} \label{eq-ellip-contour-xw}
z_i = \varphi(\alpha_j), \quad \omega_i = {1 \over N}{\diff \varphi (\alpha_j) \over \diff \alpha }, \quad \alpha_j = {2\pi  (j+1/2) \over N}, \quad j=0, \cdots, N-1.
\end{equation}
The sampling points for SS-RI are chosen as the Chebyshev points of the first kind defined in $\mathcal {I}$, which can be computed from \eqref{eq-Chebyshev-pts}. The corresponding weights are computed from \eqref{eq-Chebyshev-wts}. $N=200$ is used for both SS-RI and SS-CI.

\begin{figure}[hbt]
\centering
\epsfig{figure=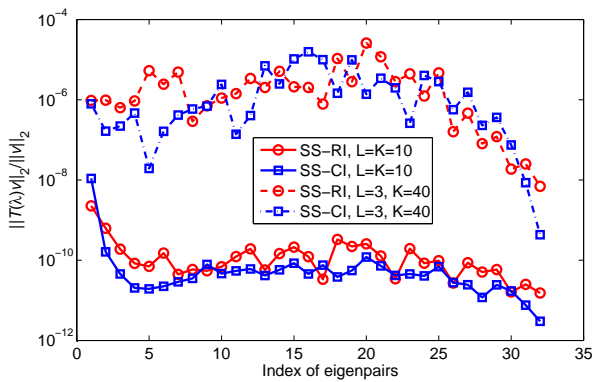,width=0.6\textwidth}
\caption{Accuracy of the SS-CI and SS-RI algorithms with different values of $L$ and $K$ in Example \ref{ex-string}}
\label{fig-cia-ria}
\end{figure}

First, we show that the two algorithms could reach the same accuracy for the same values of $N$, $K$ and $L$. To this end, we run the algorithms for two cases, `$L=3,\,K=40$' and `$L=K=10$', and
compute the relative residual $||T(\lambda)v||_2/||v||_2$ of each eigenpair $(v, \lambda)$, as demonstrated in Figure \ref{fig-cia-ria}.
When illustrate the residuals of eigenpairs, the eigenvalues are first sorted in ascending order according to their absolute values, the residual of each eigenpair is then plotted against its index in the sorted list. While having the same level of accuracy, SS-CI took roughly twice as much CPU time and memory as SS-RI due to the complex operations. For example, SS-CI used 0.92 second to compute all $T(z_i)^{-1} U$, while SS-RI used 0.49 second.

Figure 2 also shows that $L$ and $K$ have different influences on the accuracy and the algorithms would be unstable for large $K$.
When $L=K=10$, the accuracy of most eigenpairs reaches $10^{-10}$, but when $L=3$ and $K=40$ the accuracy of both algorithms becomes deteriorated, although the rank condition $K\cdot L \ge \bar n_{\mathcal {C}}$ is still satisfied. In addition, we also tried the case $L=2$ and $K=60$, the results of both algorithms were totally wrong! On the other hand, when we ran the case $L=40$ and $K=3$, the results were as good as the case $L=K=10$. In fact, $L$ is the number of vectors used to probe the resolvent $T(z)^{-1}$. It is easy to imagine that the larger the $L$, the more information of $T(z)^{-1}$ can be acquired and thus the better the results. In case of $L < \bar n_{\mathcal {C}}$, moments up to degree $K$ have to be used to guarantee the rank condition $K\cdot L \ge \bar n_{\mathcal {C}}$ via \eqref{eq-ranks}. However, since the two matrices in \eqref{eq-ranks} are Vandermonde-like matrices, they could become rank-deficient
for large $K$, as a result the rank condition \eqref{eq-ranks} would never be satisfied with good accuracy.

\end{example}

From this numerical experiment, we have shown that the RIA provides a possibility to save computational costs and memory usage for real NEPs. However, since the efficiency of the SS-RI algorithm is limited, in the following sections, we use the Rayleigh-Ritz-type algorithm \ssrrs{}
to show that the RIA allows us to enhance the solution accuracy by adjusting the positions of the sampling points.

\subsection{Superiorities of \ssrrs{} over SS-RI and \ssrrm{}} \label{S-S-nin-ssrim}

Example \ref{ex-string} is used again to study the stability and accuracy of the \ssrrs{} algorithm. Specifically,
a comparison with SS-RI will be used to show that \ssrrs{} can circumvent the instability and inaccuracy caused by using a small $L$, and a comparison with \ssrrm{} will be used to show that the sampling scheme generates better eigenspaces than the moment scheme.

\begin{figure}
\centering
\begin{subfigure}[t]{.5\textwidth}
  \centering
  \epsfig{figure=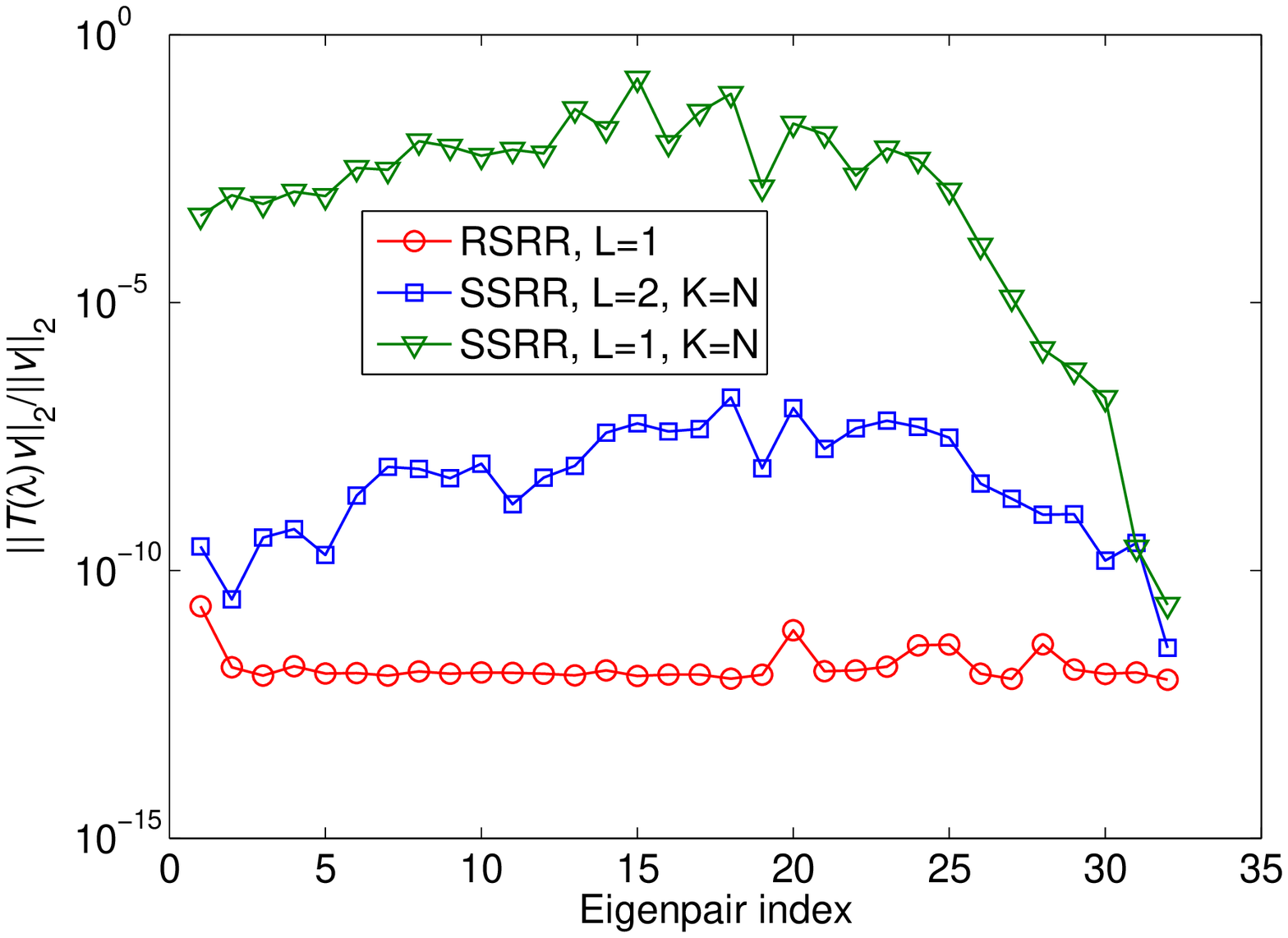,width=0.975\textwidth}
  \caption{Residuals of eigenpairs}
  \label{fig-string-rrerr}
\end{subfigure}%
\begin{subfigure}[t]{.5\textwidth}
  \centering
  \epsfig{figure=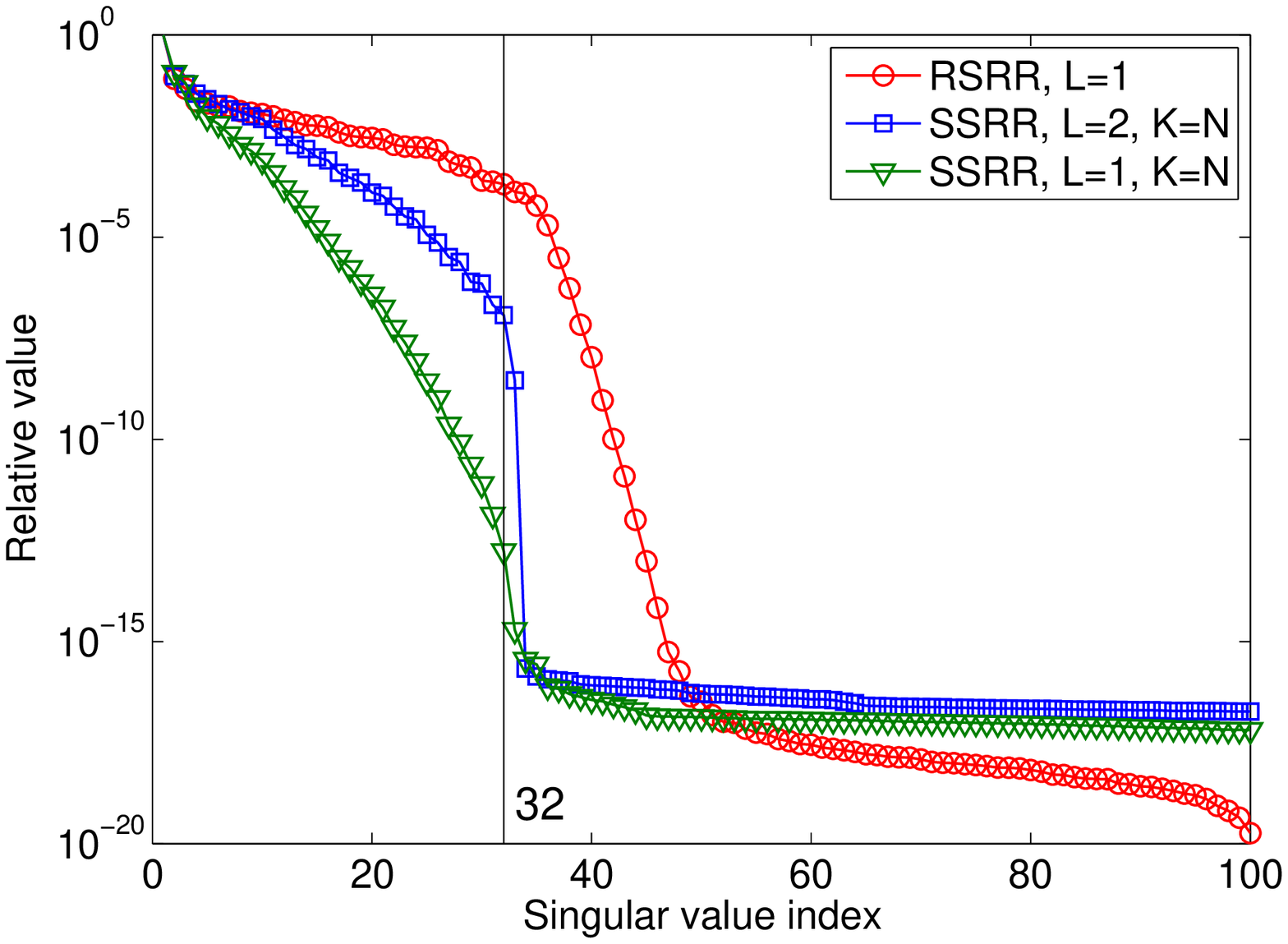,width=0.975\textwidth}
  \caption{Singular values of $S$ and $M$}
  \label{fig-string-rrsvs}
\end{subfigure}\\
\caption{Performance of \ssrrs{} and \ssrrm{} for Example \ref{ex-string}, $N=100$.}
\label{fig-string-rr}
\end{figure}

In both \ssrrs{} and \ssrrm{}, the sampling points $z_i$ are set to be Chebyshev points of the first kind in $\mathcal {I}=[3,10000]$ with $N=100$.
This corresponds to using RIA, thus only real operations are performed.
The barycentric weights $\omega_i$ of \ssrrm{} are computed from \eqref{eq-Chebyshev-wts}, and $K=N$ is used to achieve the largest size of the subspace.
Since all the eigenvalues are simple, $L=1$ is used in \ssrrs{}.
But for \ssrrm{} two cases with $L=1$ and $L=2$ are tested, respectively. The projected NEPs are solved using the SS-FULL algorithm with parameters $N_Q = 500$ and $K_Q = 8$, and the contour $\mathcal {C}_Q$ is an ellipse defined by \eqref{eq-ellip-contour} with $\gamma = 5001.5$, $a = 4998.5$, $b = 0.1a$. Note that, to avoid confusion in notations, a subscript $Q$ has been added to all the parameters of the SS-FULL algorithm that have the same meaning in the \ssrrs{} and \ssrrm{} algorithms. Figure \ref{fig-string-rr} shows the residuals of the computed eigenpairs and the singular values of matrices $S$ and $M$ in \ssrrs{} and \ssrrm{}, respectively.

Figure \ref{fig-string-rrerr} has two implications. First, by using $L=1, N=100$ \ssrrs{} can achieve an accuracy about $10^{-11}$ that is comparable to the accuracy of SS-RI with $L=10, N=200$ (see Figure \ref{fig-cia-ria}). However, the number of linear systems solved in \ssrrs{} and SS-RI are remarkable different: the former is 100, while the latter is 2000! Moreover,
when the value of $L$ is reduced to 3, SS-RI can only achieve an accuracy about $10^{-6}$ no matter how large $K$ is, and further reducing the value of $L$ leads to completely wrong results. This comparison shows that the Rayleigh-Ritz-type algorithms have much better stability and accuracy than the SS-type algorithms.

The second implication of Figure \ref{fig-string-rrerr} is about the quality of the eigenspaces constructed by the sampling scheme and the moment scheme. Since $K=N$, the eigenspaces obtained by the moment scheme reach the largest size for each $L$.
However, the residuals of the most eigenpairs of \ssrrm{} are much larger than those of \ssrrs{}.
This phenomenon can be explained by the behavior of the singular values of $S$ and $M$ in Figure \ref{fig-string-rrsvs}, where the singular values are scaled so that the first singular value $\sigma_1$ is equal to $1$.
The last 100 singular values of $M$ in case of $L=2$ have very slow decay in the interval $[10^{-17}, 10^{-18}]$ and thus are truncated.
The vertical line corresponds to the index $\bar n_\mathcal {C}$.
One can see that the singular values of $M$ drop much faster than those of $S$, due to the bad conditioning of $M$ caused by the high order moments.
In particular, the singular values of $M$ with $L=1$ drops to around $10^{-13}$ at the $\bar n_\mathcal {C}$th singular value, indicating that the rank condition \eqref{eq-rankcond-M} cannot be satisfied with a good accuracy.
Increasing $L$ improves the conditioning of $M$ so that the first $\bar n_\mathcal {C}$ singular values do not decay too fast. This leads to an improved accuracy to the eigenpairs.
In comparison, the rank property of $S$ is much better, and an obvious decay is observed only after the $\bar n_\mathcal {C}$th singular value, meaning that $S$ can fulfill the rank condition \eqref{eq-rank-NL} more securely. These results confirm the statements about the rank property of $S$ in Section \ref{S-S-samp}.

\subsection{Effects of $N$ and $L$ on the accuracy of \ssrrs{}} \label{S-S-nin-NL}

The effects of $N$ and $L$ on the computational costs of \ssrrs{} have been discussed in Section \ref{S-S-algorithm}. Here the effects of the two parameters on the accuracy of \ssrrs{} are studied using Example \ref{ex-string}. The computation is carried out for four cases with $L = 1,\, 2,\, 4$ and $8$, respectively, and in each case $N$ also increases several times.
Other setup of \ssrrs{} is the same as in Section \ref{S-S-nin-ssrim}.
For each pair of $(N,L)$, the largest relative residual of the computed eigenpairs is recorded. To compute the largest residual, the same computation is repeated 5 times and the final result is set to be the average of the 5 largest residuals.

\begin{figure}[htb]
\centering
  \epsfig{figure=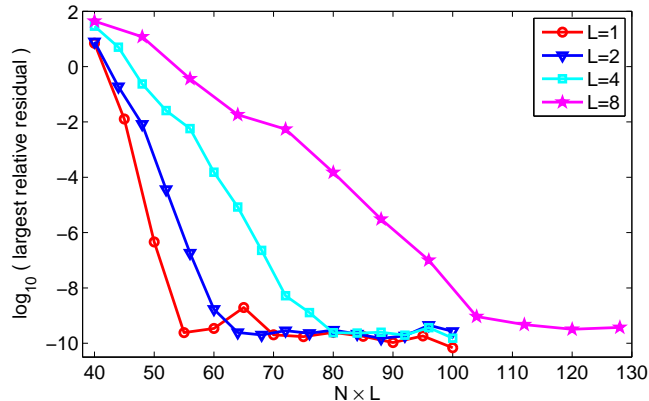,width=0.65\textwidth}
\caption{Effects of $N$ and $L$ on the accuracy of the \ssrrs{} algorithm.}
\label{fig-string-NLerr}
\end{figure}

Figure \ref{fig-string-NLerr} shows the largest residuals versus the product $N\cdot L$ for the four values of $L$. Two conclusions can be drawn.
First, for a given $L$, the error decays exponentially with the increase of $N$, and thus $N\cdot L$, till the smallest error around $10^{-10}$ is reached.
This is attributed to the fact that the weighted summation of $R_\mathcal {C}(z_i)$ \eqref{eq-filterR} converges exponentially to zero with the increase of $N$.
Second, compared with $L$, $N$ has a more remarkable influence on the accuracy of \ssrrs{}, since for a given $N\cdot L$, a larger $L$ but a smaller $N$ always corresponds to a much larger residual than a smaller $L$ but a larger $N$.

\subsection{Accuracy enhancement for \ssrrs{}} \label{S-S-accuracy-gun}

Here the ``gun'' problem of the NLEVP collection \cite{BHM13} is considered with
two purposes: (1) further compare the accuracy of \ssrrs{} and \ssrrm{} using a more general NEP with complex spectrum; (2) show that the RIA allows us to improve the accuracy of \ssrrs{} by putting sampling points close to the eigenvalues, which is not easily conceived for the CIA.

\begin{figure}[hbt]
  \centering
  \epsfig{figure=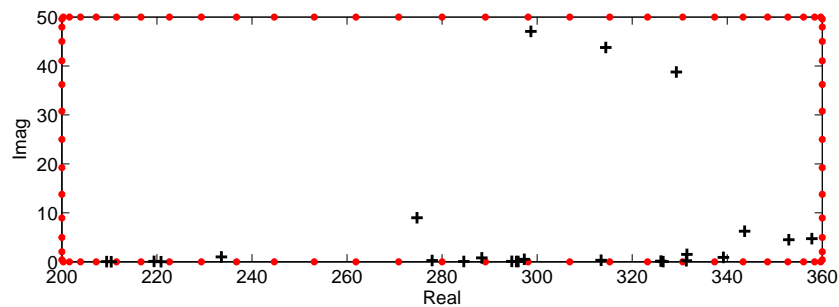,width=0.9\textwidth}
\caption{The contour (the solid line), sampling points ({\color{red}$\bullet $}) and computed eigenvalues ($\scriptstyle\pmb{+}$) by \ssrrs{} in Example \ref{ex-gun}.}
\label{fig-gun-evs}
\end{figure}

\begin{example}[\textit{Radio-frequency gun cavity} \cite{BHM13}] \label{ex-gun}
This is a large-scale problem that models a radio-frequency gun cavity, and the matrix $T(z)$ is given by \eqref{eq-nep-gun} with $J=2$,
$$
T(\lambda) = K_{\rm s}- \lambda^2 M + \im \sqrt{\lambda^2-\kappa_1^2} \,W_1 + \im \sqrt{\lambda^2-\kappa_2^2} \,W_2.
$$
The matrices $K_{\rm s}$, $M$, $W_1$ and $W_2$ are of the order $n=9956$. The cutoff values are $\kappa_1 = 0$ and $\kappa_2 = 108.8774$. We seek the 25 eigenvalues lying inside the rectangular contour $\mathcal {C}$ whose lower-left and upper-right vertices are $200 + 0 \im$ and $360+50\im$, respectively.
Figure \ref{fig-gun-evs} shows the contour and the distribution of the $25$ eigenvalues. Some of the eigenvalues are closely clustered.
This problem has been extensively used to test new algorithms for solving large-scale NEPs; see, e.g., \cite{EC13,guttel2014nleigs}. For this specific problem, we measure the convergence of an approximate eigenpair $(v, \lambda)$ by the following relative residual norm
defined in \cite{liao2010nonlinear},
$$
E(v, \lambda)= {||T(\lambda)v||_2/||v||_2 \over \|K_{\rm s}\|_1 + |\lambda^2| \|M\|_1 + \sqrt{\lambda^2-\kappa_1^2} \|W_1\|_1 + \sqrt{\lambda^2-\kappa_2^2} \|W_2 \|_1}.
$$

\end{example}

First, the performance of \ssrrs{} and \ssrrm{} is further compared.
The sampling points $z_i$ are set as the nodes of the Gauss-Legendre quadrature rules on the four sides of the rectangular $\mathcal {C}$. The 27-point rule is used on the two longer sides, and the 13-point rule is used on the two shorter sides, so totally, $N=80$. The moments of \ssrrm{} are computed by the composite Gauss quadrature rule, and $K=N$ is used. $L=2$ is used for both algorithms. The projected NEPs are solved by SS-FULL on the contour $\mathcal {C}$. The contour integrals are computed by using composite Gauss quadrature with $N_Q=500$ points, and $K_Q =2$ is used.

\begin{figure}
\centering
\begin{subfigure}[t]{.5\textwidth}
  \centering
  \epsfig{figure=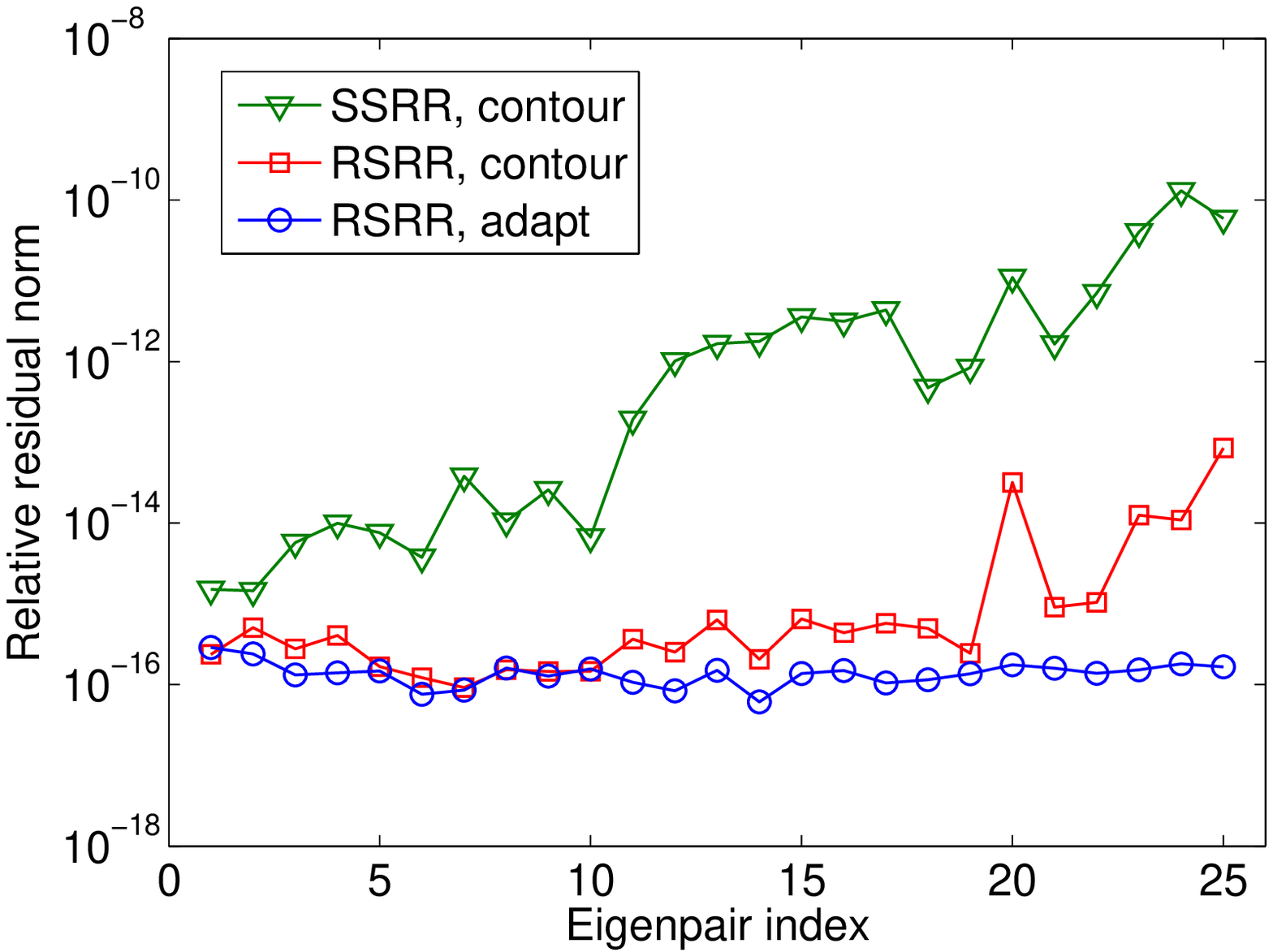,width=0.975\textwidth}
  \caption{Residuals of eigenpairs}
  \label{fig-gun-rrerr}
\end{subfigure}%
\begin{subfigure}[t]{.5\textwidth}
  \centering
  \epsfig{figure=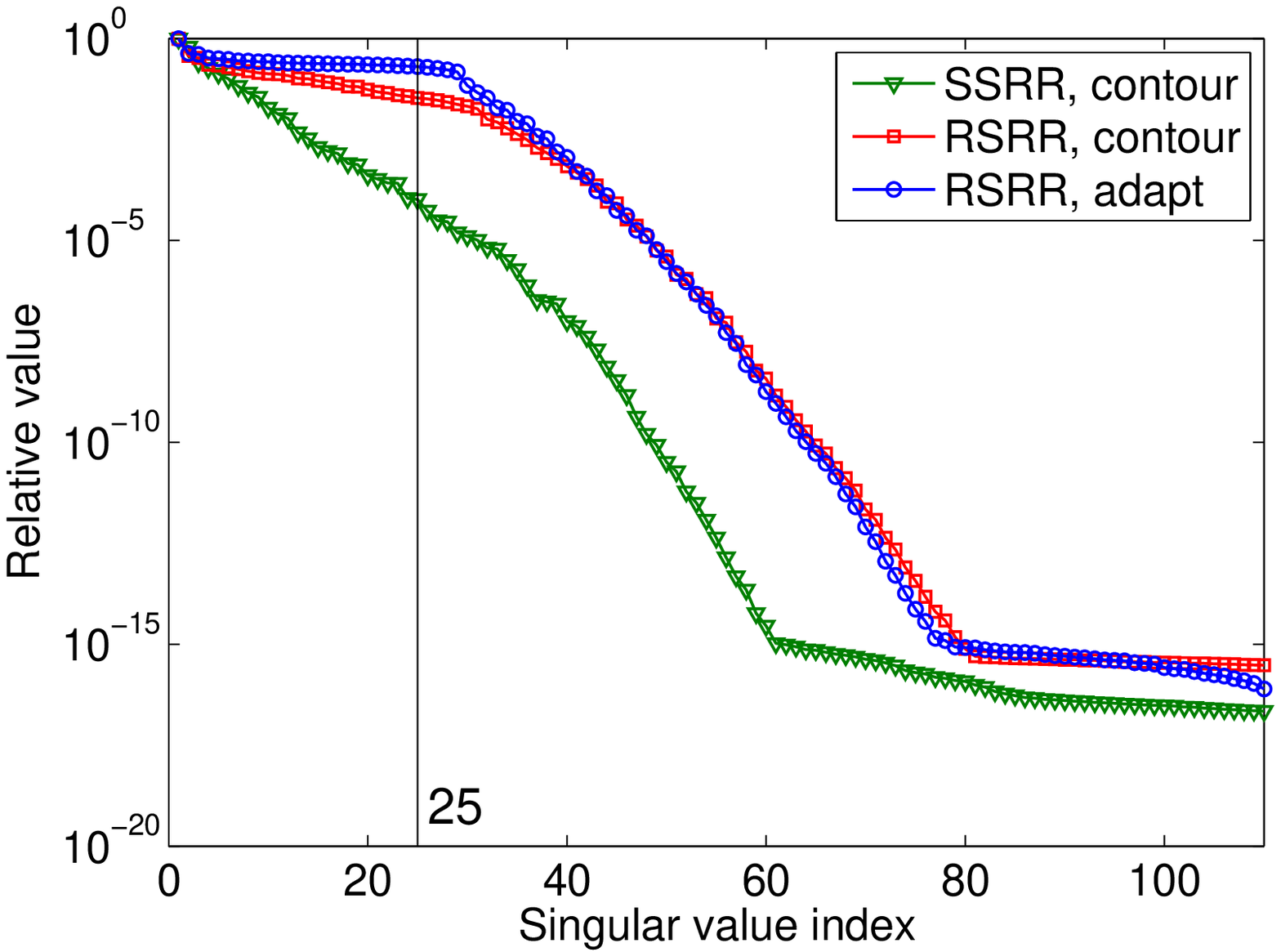,width=0.975\textwidth}
  \caption{Singular values of $S$ and $M$}
  \label{fig-gun-rrsvs}
\end{subfigure}\\
\caption{Numerical results of Example \ref{ex-gun}.}
\label{fig-gun-rr}
\end{figure}

Figure \ref{fig-gun-rr} shows the relative residual norms $E(v, \lambda)$ of the computed eigenpairs and the singular values of the matrices $S$ and $M$. The results of \ssrrs{} and \ssrrm{} are indicated by ``\ssrrs{}, contour'' and ``\ssrrm{}, contour'', respectively. The residuals of \ssrrm{} is typically more than three orders larger than those of \ssrrs{}, due to the rank-deficiency of $M$ caused by using the high order moments. This can be verified by Figure \ref{fig-gun-rrsvs}, where one can clearly see that the first $25$ singular values of $M$ decay much faster than those of $S$.

Then, we show the influence of the sampling point distribution on the accuracy of \ssrrs{}. We use a two-stage solution strategy. In the first stage, an approximation to the 25 eigenvalues is computed by \ssrrs{} using only 30 sampling points on the contour (10 points in $x$-direction and $5$ points in $y$-direction). This is a very rough approximation, and the largest residual of the 25 eigenpairs is of the order $10^{-5}$. In the second stage, the accuracy of the eigenpairs are improved by assigning the approximate eigenvalues obtained at the first stage as new sampling points.
Hence, after refinement there are totally $N=50$ ($=30+25$) sampling points.
The NEP is then solved by \ssrrs{} using these 50 sampling points, and the results are labeled by ``\ssrrs{}, adapt'' in Figure \ref{fig-gun-rr}. Comparing with the case ``\ssrrs{}, contour'' using 80 sampling points in the first test, the residuals of ``\ssrrs{}, adapt'' are typically lower and more uniformly distributed, which is consistent with the good behavior of the singular values in Figure \ref{fig-gun-rrsvs}.
This experiment shows that in order to construct better eigenspaces, it is better to put the sampling points close to the eigenvalues. The RIA offers such a possibility to distribute the sampling points. Note that when the same probing matrix $U$ is used in both the two stages, the results of $T(z)^{-1}U$ for the 30 old sampling points can be reused and thus $T(z)^{-1}U$ only has to be computed for the 25 new sampling points in the second stage.

\subsection{Comparison of \ssrrs{} with a state-of-the-art method} \label{S-S-NLEIGS-gun}

It is also meaningful to compare the performance of \ssrrs{} with some other state-of-the-art nonlinear eigensolvers. Here we choose to compare with a linearization method, NLEIGS, recently proposed in \cite{guttel2014nleigs}, because of its good performance and because its code is publicly available.

In this code, the 21 eigenvalues on the upper
half-disk centered at 223.6 with radius 111.8 are computed using four different variants of NLEIGS. We compare with the two most efficient variants: ``Variant R2'' and ``Variant S'', and the results are shown in Table \ref{tab-gun-comparison}.
By using the two-stage solution strategy in Section \ref{S-S-accuracy-gun}, \ssrrs{} computes 25 eigenpairs using 42 seconds, while ``Variant R2'' computes 21 eigenvalues using 36 seconds. Therefore, the average CPU time for one eigenpair is almost the same for these two algorithms. However, the largest residual norm of \ssrrs{} is more than two orders lower than ``Variant R2''.
``Variant S'' is the fastest one, but on the other hand its residual norm is the biggest.
In addition, the memory usage of the NLEIGS method is much more than \ssrrs{} due to the storage of the LU factorizations of the matrices $T(z_i)$ \cite{guttel2014nleigs}.

\begin{table}[!h]
\begin{center}
\caption{Results of \ssrrs{} and ``Variant R2'' of the NLEIGS method for the ``gun'' example \ref{ex-gun}.}\label{tab-gun-comparison}
\vspace{-.8\baselineskip} {\small
\begin{tabular*}{0.99\textwidth}{@{\extracolsep{\fill}}cccccc@{}}\toprule
      Algorithms           &  \#$\lambda$s &  Max. $E(v, \lambda)$  & CPU time & Time per $\lambda$ & Memory usage     \\ \hline
      ``\ssrrs{}, adapt''  &  25           &  $3 \times 10^{-16}$   & 42 s     & 1.68               & $\thicksim$140 MB\\
      ``Variant R2''       &  21           &  $1 \times 10^{-13}$   & 36 s     & 1.71               & $\thicksim$450 MB\\
      ``Variant S''        &  21           &  $1 \times 10^{-11}$   & 13 s     & 0.62               & $\thicksim$435 MB\\
\bottomrule
\end{tabular*}}
\end{center}
\end{table}

From the above results, we see that the overall performance of \ssrrs{} is better than the ``Variant R2'' of the NLEIGS method. ``Variant S'' is the fastest approach for this specific example, but as shown in \cite{guttel2014nleigs}, the performance of ``Variant S'' is sensitive to the number and distribution of the shifts in the rational Krylov method, and an example is given in \cite{guttel2014nleigs} to show that the ``Variant S'' can be inferior to the ``Variant R2'' in some cases. Lastly,
we stress that the \ssrrs{} algorithm is easier to implement and parallelize than the NLEIGS method.

\section{Engineering applications and validation} \label{S-ne-realapps}

Two more realistic examples are provided to demonstrate the performance of the newly developed \ssrrs{} algorithm in dealing with engineering problems. The first example comes from the FEM analysis of a viscoelastically damped payload attach fitting structure, which has about 1 million degrees of freedom (DOFs). The second example is the BEM modal analysis of the sound field in a car cabin, which has around 60000 DOFs.

\begin{figure}[htb]
\setlength{\unitlength}{0.0108889\textwidth}
\centering
\begin{overpic}[width=0.49\textwidth]%
{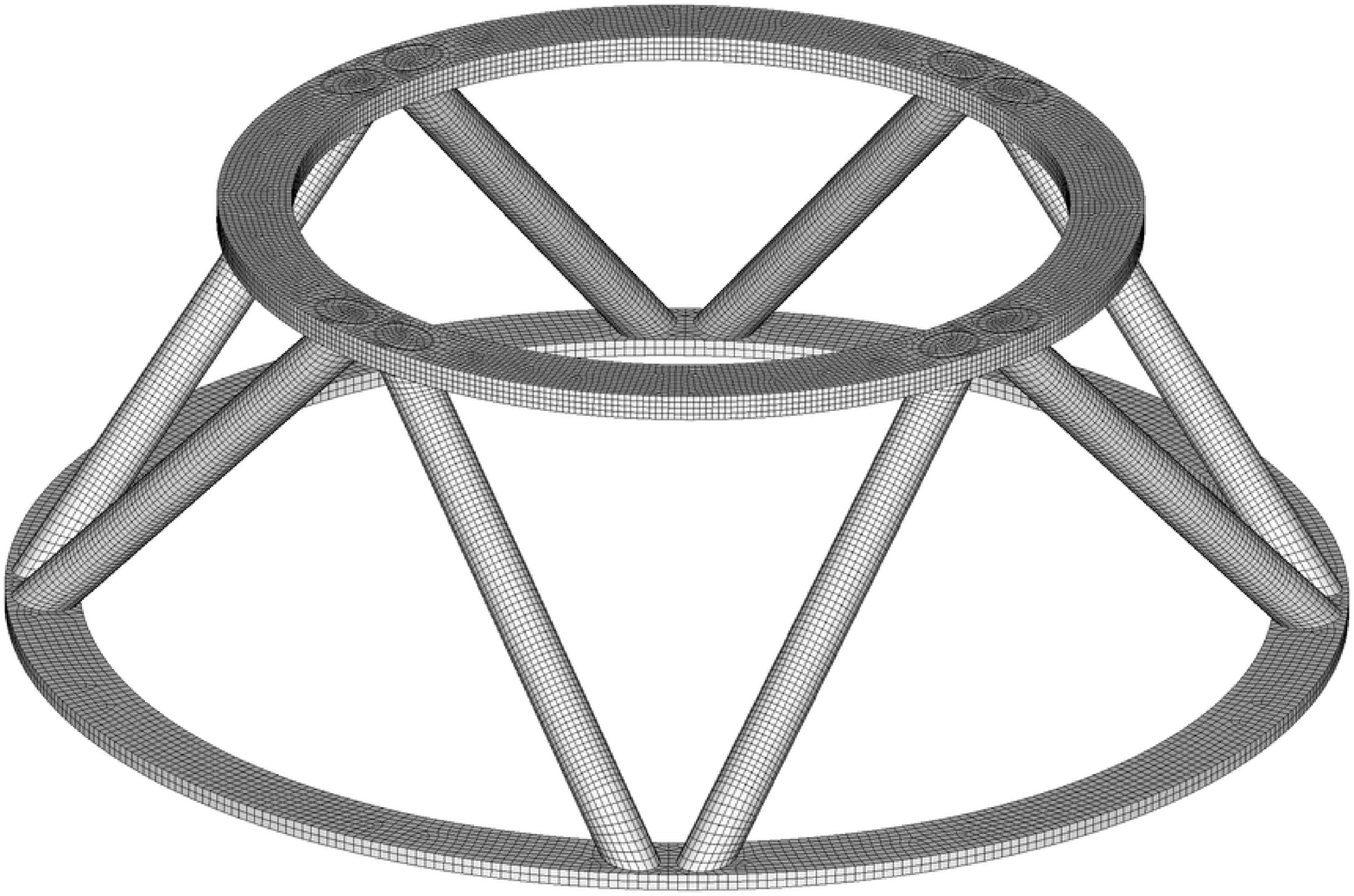}

\begin{picture}(1,1)

\put(21,10.5){\vector(1,0){3}}
\put(21,10.5){\vector(0,1){3}}
\put(21,10.5){\vector(-1,-1){2.6}}
\put(23.5,11.5){{\footnotesize $y$}}
\put(20,14){\footnotesize $z$}
\put(18,9.5){\footnotesize $x$}

\put(10.5,7){\vector(-1,-1){2.6}}
\put(4,8){{\tiny Lower Flange}}
\put(4,28){\vector(3,-2){2.6}}
\put(-2,28.8){{\tiny Upper Flange}}
\put(40,22){\vector(-1,-1){2.6}}
\put(37,24.5){{\tiny Damping}}
\put(37,23){{\tiny Cylinder}}
\end{picture}
\end{overpic}\hspace{0.5cm}
\epsfig{figure=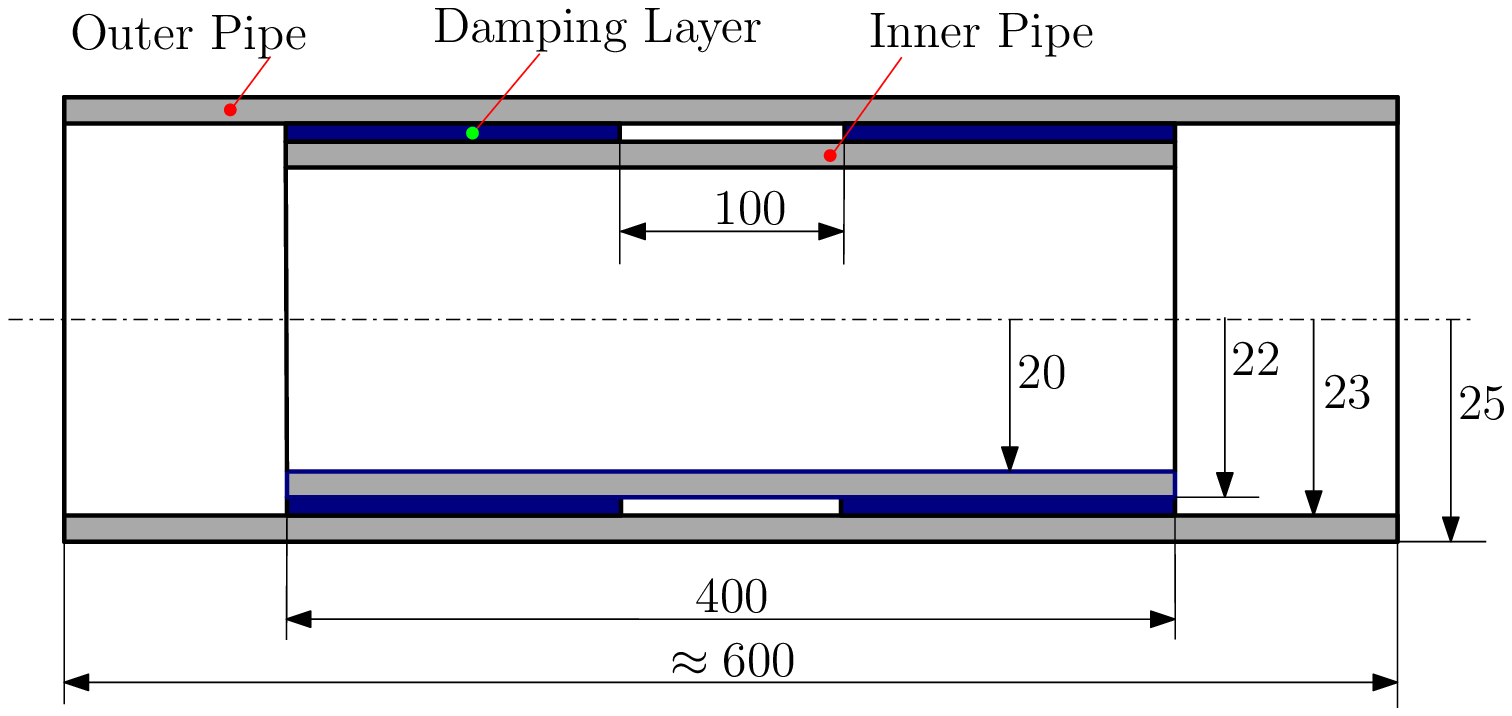,width=0.47\textwidth}
\caption{The payload attach fitting model (Unit: mm). The model consists of an upper flange, a lower flange and eight damping cylinders. The dimensions of the upper flange are: inner radius 365, outer radius 435, thickness 24. The dimensions of the lower flange are: inner radius 560, outer radius 630, thickness 16. The structure and dimensions of the damping cylinders are shown in the left diagram.}
\label{fig-paf-model}
\end{figure}

\begin{example}[\textit{Payload attach fitting}] \label{ex-paf}

The payload attach fitting structures often serve as isolators
between the satellite and the launch vehicle to reduce the satellite vibration caused by the launch-induced dynamic loads. Here the model illustrated in Figure \ref{fig-paf-model} is considered.
The structural material is Aluminum with elastic modulus $70\,\rm{GPa}$, density $2770\, \rm{kg}/\rm{m}^3$ and Poisson's ratio $0.3$. The eight damping cylinders contain viscoelastic damping material ZN1 polymer which is modeled by the Biot model. The density and Poisson's ratio of the damping material are $970\, \rm{kg}/\rm{m}^3$ and $0.49$, respectively. The structure is fixed at the lower surface of the lower flange and all the other surfaces are free of traction.

This problems is solved by using the FEM software \textsc{Ansys}\textregistered{}.
The model is discretized by using the SOLID186 element, and the total number of DOFs is 1005648. The expression of this NEP is of the form \eqref{eq-REP-damping} with $T(z) = z^2 M + zG(z)K_{\rm v} + K_{\rm s}$,
$zG(z) = G_{\infty} \left(1 + \sum_{k=1}^4 {a_k z \over z + b_k}\right)$, $G_{\infty} = 362750\,\rm{Pa}$, $a_1 = 0.762063$, $a_2 = 1.814626$, $a_3 = 84.93828$, $a_4 = 4.869723$, $b_1 = 53.72964$, $b_2 = 504.5871$, $b_3 = 29695.64$ and $b_4 = 2478.43$. The natural frequencies below 800Hz are sought. A rough modal analysis for the undamped structure shows that the lowest natural frequency is around 204 Hz, hence the frequency interval for \ssrrs{} is set as $[200, 800]$Hz. The computation was carried out on a Server with eight 8-core Intel Xeon E7-8837 (2.67GHz) processors and 256 GB RAM.

\end{example}

 In order to take the full advantage of the RIA and to keep the total number of sampling points reasonably small, this problem is solved using the two-stage strategy in Section \ref{S-S-nin-gun}. In the first stage, the eigenvalues are estimated by using only $20$ sampling points (Chebyshev points of the first kind) in the interval $[200, 800]$Hz. In the second stage, eigenvalues computed before are set as new sampling points, and then the problem is solved again using both the old and new sampling points.

In all the \ssrrs{} solutions, $L = 2$ is used. The linear systems involved in computing $T(z)^{-1}U$ are solved by using the default sparse direct solver (the frontal solver with $LDL^T$ factorization) in \textsc{Ansys}\textregistered{} without any special treatment.
The reduced NEPs are solved by using SS-FULL on the ellipse contour \eqref{eq-ellip-contour} with $\gamma = 500, \, a = 300$ and $b = 0.5a$. The contour integrals are computed by using the trapezoidal rule \eqref{eq-ellip-contour-xw} with $N_Q=500$ points, and $K_Q =2$ is used.
The entire solution process takes about 2 hours.

Table \ref{tab-paf-evs} shows the computed eigenvalues and residuals in both the two stages. This problem has 25 eigenvalues, as listed in the third column, but only 24 eigenpairs are obtained by using the 20 sampling points in the first stage, including 5 2-fold eigenvalues. The largest relative residual of the 24 eigenpairs reaches 4.415. In the second stage, the 19 mutually different eigenvalues are assigned as new sampling points and thus $N = 39$ sampling points are used. The \ssrrs{} algorithm correctly obtains the 25 eigenvalues with a remarkably improved accuracy: the largest residual is about $2.663\times 10^{-6}$, six orders lower than that of the first stage solution.

The same payload attach fitting model has also been solved by using the CIA based \ssrrs{} in \cite{xiao2015contour}, where 14 eigenpairs in the smaller interval $[160, 640]$ are computed using $N=100$ sampling points on an elliptical contour. The largest residual of the computed eigenpairs is about $10^{-5}$. Whereas, here 25 eigenpairs are solved using only $N=39$ sampling points, but better results have been obtained.
This example confirms again that based on the newly developed RIA, it is possible to considerably improve the accuracy of the \ssrrs{} algorithm by a judicious selection of the sampling points.

\begin{table}[!h]
\begin{center}
\caption{Computed eigenvalues and residuals of the eigenpairs in the two solution stages for the payload attach fitting model}\label{tab-paf-evs}
\vspace{-.8\baselineskip} {\small
\begin{tabular*}{0.99\textwidth}{@{\extracolsep{\fill}}lllll@{}}\toprule
& \multicolumn{2}{c}{Computed eigenvalues} & \multicolumn{2}{c}{$||T(\lambda)v||_2/||v||_2$ }  \\
\cline{2-3}  \cline{4-5}
$k$         &      First stage   & Second stage        &  First stage   & Second stage  \\
\hline
1  & 206.877054534+0.993274962i	& 206.877054436+0.993275048i   & 4.570E-6	& 2.727E-8  \\
2  & 228.061075395+1.788371341i	& 228.061075870+1.788370624i   & 9.162E-6	& 1.064E-6  \\
3  & 228.061075854+1.788370607i	& 228.061075876+1.788370623i   & 2.350E-7	& 1.749E-6 \\
4  & 251.945403991+0.164908920i	& 251.945404108+0.164908807i   & 5.258E-5	& 3.962E-8 \\
5  & 289.734769613+2.908555918i	& 289.734769656+2.908555874i   & 3.644E-7	& 2.915E-8 \\
6  & 345.704771998+12.245151119i& 345.704779972+12.245152620i  & 2.407E-5   & 1.257E-6 \\
7  & 345.704780142+12.245152448i& 345.704779973+12.245152607i  & 6.828E-7	& 2.124E-6 \\
8  & 405.560427505+8.454299109i	& 405.560414500+8.454230992i   & 3.410E-4	& 4.491E-8\\
9  & 431.312424354+12.255687576i& 431.312295649+12.255477928i  & 3.107E-3	& 4.322E-8\\
10 & 551.306138347+6.256446835i	& 551.306173678+6.256449037i   & 5.499E-4	& 7.713E-8\\
11 & 551.307091951+6.258569739i	& 551.306173684+6.256449057i   & 1.468E-2	& 3.885E-7\\
12 & 573.462580570+2.975001736i	& 573.462544813+2.975002275i   & 4.061E-4	& 6.476E-8\\
13 & 607.309221542+15.520360915i& 607.297242169+15.478799600i  & 4.439E-1	& 5.275E-8\\
14 & 631.260975030+17.832595421i& 630.172433959+17.915518908i  & 4.145E-0	& 5.741E-7\\
15 & 640.180739896+7.657545670i	& 640.180790794+7.657485207i   & 1.016E-3	& 8.588E-7  \\
16 & 640.227479171+7.576989109i	& 640.180790805+7.657485206i   & 4.760E-1	& 1.156E-6  \\
17 & 688.454146650+2.982724282i	& 688.454340293+2.982729517i   & 4.374E-3	& 1.406E-6 \\
18 & 689.334173861+2.566616866i	& 688.454340307+2.982729542i   & 6.231E-1	& 3.916E-7 \\
19 &           ------           & 690.858208687+1.821700043i   & ------     & 2.663E-6 \\
20 & 704.728697215+3.157964964i	& 705.107604497+3.129254146i   & 1.484E-0	& 1.130E-7 \\
21 & 719.385017165+3.921698409i	& 719.385375150+3.923570746i   & 1.833E-2	& 5.591E-8 \\
22 & 723.421877396+3.240014095i	& 723.444287703+3.228000754i   & 2.084E-1	& 3.719E-7\\
23 & 723.459627245+3.218688142i	& 723.444287708+3.228000763i   & 1.946E-1	& 1.134E-7\\
24 & 742.077016154+4.995761542i	& 742.132722448+4.768007224i   & 1.954E-0 	& 1.141E-7\\
25 & 746.455479936+1.908005953i	& 746.542462488+1.931098174i   & 6.956E-1	& 9.495E-8\\
\bottomrule
\end{tabular*}}
\end{center}
\end{table}

\begin{example}[\textit{Car cabin}] \label{ex-car}

As the last example, \ssrrs{} is used to solve a NEP arising from the BEM acoustic modal analysis of a car cabin cavity. The BEM analysis is carried out using our in-house code that implements the fast BEM in \cite{CWX15}. The boundary of the cabin cavity is partitioned into 9850 triangular quadratic elements, so the BEM model consists of 59100 DOFs; see Figure \ref{fig-car-mesh} for the model shape and mesh. The entire boundary is assumed to be rigid. Such problems are challenging for the current BEM eigensolvers; see Section \ref{S-S-typicalNEPs} and \cite{xiao2015contour}.

Eigenvalues in the real interval $[40,500]$Hz are sought. The sampling points are chosen as the Chebyshev points in the interval; $L = 2$, $N = 200$ are used. The reduced NEP is solved by SS-FULL, with the sampling points being the $N_Q = 1000$ Chebyshev points in the interval, and $K_Q = 2$.
The linear systems $T(z_i)^{-1} U,\, i = 0, \cdots, N-1$ are solved by the GMRES solver with ILU preconditioner. Both the accuracy of the fast BEM and the convergence tolerance of the GMRES solver are set to be $10^{-6}$. The matrix $T(z)$ is interpolated by using Chebyshev polynomials of degree $d = 40$; see \cite{xiao2015contour} for the details. The computation was performed on a personal computer with an Intel$^\circledR$
Core$^\mathrm{TM}$ i3-2100 (3.10 GHz) CPU and 16 GB RAM. 

\begin{figure}[htb]
\centering
  \epsfig{figure=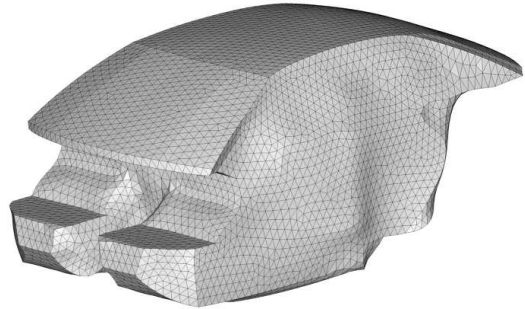,width=0.5\textwidth}
\caption{BEM mesh for the car cabin model in Example \ref{ex-car}}
\label{fig-car-mesh}
\end{figure}

\begin{figure}[htb]
\centering
\epsfig{figure=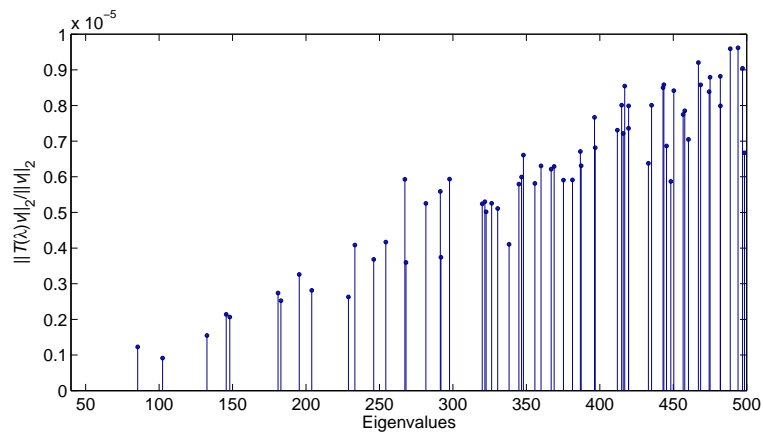,width=0.8\textwidth}
\caption{The computed eigenvalues and residuals of eigenpairs for Example \ref{ex-car}.}
\label{fig-car-evs}
\end{figure}

\ssrrs{} returns 64 eigenvalues, which are indicated by the horizontal axis in Figure \ref{fig-car-evs}. The eigenvalues are more intensively distributed in high frequency band.
The vertical axis in Figure \ref{fig-car-evs} represents the residuals of the corresponding eigenpairs.
All the residuals are below $1.0\times 10^{-5}$. The computations for this problem takes about 24 hours, with 17 hours spending on the solution of the $N\cdot L = 400$ linear systems, and around 7 hours on the computation of the Chebyshev interpolation of the reduced matrix $T_Q (z)$. Other operations, such as the computation of the truncated SVD of $S$ and running SS-FULL, only takes several minutes. The total memory usage is about 2.2 GB.

\end{example}

\section{Conclusions} \label{S-conclusions}

The paper has been devoted to the development of robust numerical methods for solving large-scale NEPs in science and engineering. In particular, the NEPs in the FEM and BEM are considered and solved. The outcomes lie in two aspects.

Theoretically, a framework for solving NEPs, called \emph{rational interpolation approach}, has been developed. It encloses the existing contour integral approach as a special case, and provides more possibility to select the sampling points $z_i$, at which the resolvent $T(z)^{-1}U$ is computed. For general complex NEPs, it allows us to improve the accuracy of the eigenpairs by putting the sampling points close to the eigenvalues. For real NEPs, it allows us to reduce the computational time and memory consumption by using real sampling points and performing real arithmetics.
Besides, it inherits the two salient merits of the contour integral approach, that is, it simultaneously considers all the eigenvalues within a given contour, and enjoys good parallelizability.

Numerically, a robust and accurate NEP solver, called \emph{\ssrrs{}}, has been proposed based on the rational interpolation approach. The success of \ssrrs{} is established on the newly proposed sampling scheme for generating approximate eigenspaces. Unlike the existing moment scheme, where the moments of the resolvent $T(z)^{-1}U$ are used to construct eigenspaces and thus may fail if high order moments are involved, the sampling scheme generates eigenspaces directly from the values of $T(z)^{-1}U$ at the sampling points $z_i$. The eigenspaces of the sampling scheme are larger and more reliable, leading to the higher accuracy of the \ssrrs{} algorithm. The applicability of \ssrrs{} is not limited by the structure of the system matrix $T(z)$ and the characteristics of the eigen-solutions. Moreover, \ssrrs{} can be easily implemented and parallelized.

The good robustness and accuracy of the \ssrrs{} algorithm have been demonstrated by a variety of typical examples and comparisons with other nonlinear eigensolvers. In particular, a FEM model with around 1 million DOFs and a BEM model with several tens of thousands of DOFs are supplied to show the efficiency of \ssrrs{} in solve large-scale engineering NEPs.

\section*{Acknowledgements}

JX gratefully acknowledges the financial supports from the National Science
Foundations of China under Grants 11102154 and 11472217, Fundamental Research Funds for the Central Universities and the Alexander von Humboldt Foundation (AvH) to support his fellowship research at the Chair of Structural Mechanics, University of Siegen, Germany.

\bibliographystyle{unsrt}
\bibliography{nep_ii_bib}

\begin{thebibliography}{10}

\bibitem{TM01}
Fran{\c{c}}oise Tisseur and Karl Meerbergen.
\newblock The quadratic eigenvalue problem.
\newblock {\em SIAM Review}, 43(2):235--286, 2001.

\bibitem{MS11}
Volker Mehrmann and Christian Schr{\"o}der.
\newblock Nonlinear eigenvalue and frequency response problems in industrial
  practice.
\newblock {\em Journal of Mathematics in Industry}, 1(1):1--18, 2011.

\bibitem{effenberger2013robust}
Cedric Effenberger.
\newblock {\em Robust solution methods for nonlinear eigenvalue problems}.
\newblock PhD thesis, {\'E}cole polytechnique f{\'e}d{\'e}rale de Lausanne,
  2013.

\bibitem{van2015rational}
Roel Van~Beeumen.
\newblock {\em Rational {Krylov} methods for nonlinear eigenvalue problems}.
\newblock PhD thesis, KU Leuven, 2015.

\bibitem{OBZ15}
TM~van Opstal, EH~van Brummelen, and GJ~van Zwieten.
\newblock A finite-element/boundary-element method for three-dimensional,
  large-displacement fluid--structure-interaction.
\newblock {\em Computer Methods in Applied Mechanics and Engineering},
  284:637--663, 2015.

\bibitem{Steinbach14K}
A~Kimeswenger, O~Steinbach, and G~Unger.
\newblock Coupled finite and boundary element methods for fluid-solid
  interaction eigenvalue problems.
\newblock {\em SIAM Journal on Numerical Analysis}, 52(5):2400--2414, 2014.

\bibitem{abdel1994safeguarded}
Mohammedi~R Abdel-Aziz.
\newblock Safeguarded use of the implicit restarted {Lanczos} technique for
  solving non-linear structural eigensystems.
\newblock {\em International Journal for Numerical Methods in Engineering},
  37(18):3117--3133, 1994.

\bibitem{dumont2007solution}
NA~Dumont.
\newblock On the solution of generalized non-linear complex-symmetric
  eigenvalue problems.
\newblock {\em International Journal for Numerical Methods in Engineering},
  71(13):1534--1568, 2007.

\bibitem{EC13}
C~Effenberger.
\newblock Robust successive computation of eigenpairs for nonlinear eigenvalue
  problems.
\newblock {\em SIAM Journal on Matrix Analysis and Applications},
  34(3):1231--1256, 2013.

\bibitem{salas2014spectral}
Pablo Salas, Luc Giraud, Yousef Saad, and St{\'e}phane Moreau.
\newblock Spectral recycling strategies for the solution of nonlinear
  eigenproblems in thermoacoustics.
\newblock {\em Preprint ys-2015-1, Dept. Computer Science and Engineering,
  University of Minnesota}, 2015.

\bibitem{lu2015pade}
Ding Lu, Xin Huang, Zhaojun Bai, and Yangfeng Su.
\newblock A {Pad{\'e}} approximate linearization algorithm for solving the
  quadratic eigenvalue problem with low-rank damping.
\newblock {\em International Journal for Numerical Methods in Engineering},
  2015.

\bibitem{sakurai2003projection}
Tetsuya Sakurai and Hiroshi Sugiura.
\newblock A projection method for generalized eigenvalue problems using
  numerical integration.
\newblock {\em Journal of Computational and Applied Mathematics},
  159(1):119--128, 2003.

\bibitem{polizzi2009density}
Eric Polizzi.
\newblock Density-matrix-based algorithm for solving eigenvalue problems.
\newblock {\em Physical Review B}, 79(11):115112, 2009.

\bibitem{AS09}
Junko Asakura, Tetsuya Sakurai, Hiroto Tadano, Tsutomu Ikegami, and Kinji
  Kimura.
\newblock A numerical method for nonlinear eigenvalue problems using contour
  integrals.
\newblock {\em JSIAM Letters}, 1(0):52--55, 2009.

\bibitem{Beyn12}
Wolf-J{\"u}rgen Beyn.
\newblock An integral method for solving nonlinear eigenvalue problems.
\newblock {\em Linear Algebra and Its Applications}, 436(10):3839--3863, 2012.

\bibitem{YS13}
Shinnosuke Yokota and Tetsuya Sakurai.
\newblock A projection method for nonlinear eigenvalue problems using contour
  integrals.
\newblock {\em JSIAM Letters}, 5(0):41--44, 2013.

\bibitem{gavin2013non}
Brendan Gavin and Eric Polizzi.
\newblock Non-linear eigensolver-based alternative to traditional {SCF}
  methods.
\newblock {\em The Journal of Chemical Physics}, 138(19):194101, 2013.

\bibitem{austin2015computing}
Anthony~P Austin and Lloyd~N Trefethen.
\newblock Computing eigenvalues of real symmetric matrices with rational
  filters in real arithmetic.
\newblock {\em SIAM Journal on Scientific Computing}, 37(3):A1365--–A1387,
  2015.

\bibitem{stefan2015zolotarev}
Stefan G{\"u}ttel, Eric Polizzi, Ping Tak~Peter Tang, and Gautier Viaud.
\newblock Zolotarev quadrature rules and load balancing for the {FEAST}
  eigensolver.
\newblock {\em SIAM Journal on Scientific Computing}, 37(4):A2100–--A2122,
  2015.

\bibitem{SAT09}
Tetsuya Sakurai, Junko Asakura, Hiroto Tadano, and Tsutomu Ikegami.
\newblock Error analysis for a matrix pencil of {Hankel} matrices with
  perturbed complex moments.
\newblock {\em JSIAM Letters}, 1(0):76--79, 2009.

\bibitem{VRK13}
Roel Van~Beeumen, Karl Meerbergen, and Wim Michiels.
\newblock A rational {Krylov} method based on {Hermite} interpolation for
  nonlinear eigenvalue problems.
\newblock {\em SIAM Journal on Scientific Computing}, 35(1):A327--A350, 2013.

\bibitem{Neu85}
A~Neumaier.
\newblock Residual inverse iteration for the nonlinear eigenvalue problem.
\newblock {\em SIAM Journal on Numerical Analysis}, 22(5):914--923, 1985.

\bibitem{Voss07}
H~Voss.
\newblock A {Jacobi--Davidson} method for nonlinear and nonsymmetric
  eigenproblems.
\newblock {\em Computers \& Structures}, 85(17):1284--1292, 2007.

\bibitem{kressner2009block}
Daniel Kressner.
\newblock A block {Newton} method for nonlinear eigenvalue problems.
\newblock {\em Numerische Mathematik}, 114(2):355--372, 2009.

\bibitem{hochstenbach2009controlling}
Michiel~E Hochstenbach and Yvan Notay.
\newblock Controlling inner iterations in the {Jacobi-Davidson} method.
\newblock {\em SIAM Journal on Matrix Analysis and Applications},
  31(2):460--477, 2009.

\bibitem{effenberger2012linearization}
Cedric Effenberger, Daniel Kressner, and Christian Engstr{\"o}m.
\newblock Linearization techniques for band structure calculations in absorbing
  photonic crystals.
\newblock {\em International Journal for Numerical Methods in Engineering},
  89(2):180--191, 2012.

\bibitem{kirkup1993solution}
Stephen~Martin Kirkup and S~Amini.
\newblock Solution of the {Helmholtz} eigenvalue problem via the boundary
  element method.
\newblock {\em International Journal for Numerical Methods in Engineering},
  36(2):321--330, 1993.

\bibitem{ARY95}
A~Ali, C~Rajakumar, and SM~Yunus.
\newblock Advances in acoustic eigenvalue analysis using boundary element
  method.
\newblock {\em Computers \& Structures}, 56(5):837--847, 1995.

\bibitem{EK12}
Cedric Effenberger and Daniel Kressner.
\newblock Chebyshev interpolation for nonlinear eigenvalue problems.
\newblock {\em BIT Numerical Mathematics}, 52(4):933--951, 2012.

\bibitem{van2015compact}
Roel Van~Beeumen, Karl Meerbergen, and Wim Michiels.
\newblock Compact rational {Krylov} methods for nonlinear eigenvalue problems.
\newblock {\em SIAM Journal on Matrix Analysis and Applications},
  36(2):820--838, 2015.

\bibitem{mackey2006structured}
D~Steven Mackey, Niloufer Mackey, Christian Mehl, and Volker Mehrmann.
\newblock Structured polynomial eigenvalue problems: {Good} vibrations from
  good linearizations.
\newblock {\em SIAM Journal on Matrix Analysis and Applications},
  28(4):1029--1051, 2006.

\bibitem{EKSU12}
Cedric Effenberger, Daniel Kressner, Olaf Steinbach, and Gerhard Unger.
\newblock Interpolation-based solution of a nonlinear eigenvalue problem in
  fluid-structure interaction.
\newblock {\em PAMM}, 12(1):633--634, 2012.

\bibitem{quraishi2014solution}
Sarosh Quraishi, Christian Schr{\"o}der, and Volker Mehrmann.
\newblock Solution of large scale parametric eigenvalue problems arising from
  brake squeal modeling.
\newblock {\em Proceedings in Applied Mathematics and Mechanics}, 2014.

\bibitem{SSI06}
Sergey~I Solov'{\"e}v.
\newblock Preconditioned iterative methods for a class of nonlinear eigenvalue
  problems.
\newblock {\em Linear Algebra and Its Applications}, 415(1):210--229, 2006.

\bibitem{conca1989existence}
Carlos Conca, J~Planchard, and M~Vanninathan.
\newblock Existence and location of eigenvalues for fluid-solid structures.
\newblock {\em Computer Methods in Applied Mechanics and Engineering},
  77(3):253--291, 1989.

\bibitem{trindade2000modeling}
MA~Trindade, A~Benjeddou, and R~Ohayon.
\newblock Modeling of frequency-dependent viscoelastic materials for
  active-passive vibration damping.
\newblock {\em Journal of Vibration and Acoustics}, 122(2):169--174, 2000.

\bibitem{liao2010nonlinear}
Ben-Shan Liao, Zhaojun Bai, Lie-Quan Lee, and Kwok Ko.
\newblock Nonlinear {Rayleigh-Ritz} iterative method for solving large scale
  nonlinear eigenvalue problems.
\newblock {\em Taiwanese Journal of Mathematics}, 14(3A):869--883, 2010.

\bibitem{nicoud2007acoustic}
Franck Nicoud, Laurent Benoit, Claude Sensiau, and Thierry Poinsot.
\newblock Acoustic modes in combustors with complex impedances and
  multidimensional active flames.
\newblock {\em AIAA journal}, 45(2):426--441, 2007.

\bibitem{CWX15}
Yanchuang Cao, Lihua Wen, Jinyou Xiao, and Yijun Liu.
\newblock A fast directional {BEM} for large-scale acoustic problems based on
  the {Burton--Miller} formulation.
\newblock {\em Engineering Analysis with Boundary Elements}, 50:47--58, 2015.

\bibitem{van2015designing}
Marc Van~Barel.
\newblock Designing rational filter functions for solving eigenvalue problems
  by contour integration.
\newblock {\em Linear Algebra and its Applications, in press}, 2015.

\bibitem{austin2014numerical}
Anthony~P Austin, Peter Kravanja, and Lloyd~N Trefethen.
\newblock Numerical algorithms based on analytic function values at roots of
  unity.
\newblock {\em SIAM Journal on Numerical Analysis}, 52(4):1795--1821, 2014.

\bibitem{saff1972extension}
EB~Saff.
\newblock An extension of {Montessus} de {Ballore's} theorem on the convergence
  of interpolating rational functions.
\newblock {\em Journal of Approximation Theory}, 6(1):63--67, 1972.

\bibitem{eugeciouglu1989fast}
{\"O}mer E{\u{g}}ecio{\u{g}}lu and {\c{C}}etin~K Ko{\c{c}}.
\newblock A fast algorithm for rational interpolation via orthogonal
  polynomials.
\newblock {\em Mathematics of Computation}, 53(187):249--264, 1989.

\bibitem{berrut2004barycentric}
Jean-Paul Berrut and Lloyd~N Trefethen.
\newblock Barycentric {Lagrange} interpolation.
\newblock {\em SIAM Review}, 46(3):501--517, 2004.

\bibitem{sakurai2013efficient}
Tetsuya Sakurai, Yasunori Futamura, and Hiroto Tadano.
\newblock Efficient parameter estimation and implementation of a contour
  integral-based eigensolver.
\newblock {\em Journal of Algorithms \& Computational Technology},
  7(3):249--269, 2013.

\bibitem{wang2014explicit}
Haiyong Wang, Daan Huybrechs, and Stefan Vandewalle.
\newblock Explicit barycentric weights for polynomial interpolation in the
  roots or extrema of classical orthogonal polynomials.
\newblock {\em Mathematics of Computation}, 83(290):2893--2914, 2014.

\bibitem{parks2006recycling}
Michael~L Parks, Eric De~Sturler, Greg Mackey, Duane~D Johnson, and Spandan
  Maiti.
\newblock Recycling krylov subspaces for sequences of linear systems.
\newblock {\em SIAM Journal on Scientific Computing}, 28(5):1651--1674, 2006.

\bibitem{hasegawa2015recovering}
Tetsuya Hasegawa, Akira Imakura, and Tetsuya Sakurai.
\newblock Recovering from accuracy deterioration in the contour integral-based
  eigensolver.
\newblock {\em JSIAM Letters}, 8(0):1--4, 2015.

\bibitem{BHM13}
Timo Betcke, Nicholas~J Higham, Volker Mehrmann, Christian Schr{\"o}der, and
  Fran{\c{c}}oise Tisseur.
\newblock {NLEVP}: A collection of nonlinear eigenvalue problems.
\newblock {\em ACM Transactions on Mathematical Software (TOMS)}, 39(2):7,
  2013.

\bibitem{guttel2014nleigs}
Stefan G\"uttel, Roel Van~Beeumen, Karl Meerbergen, and Wim Michiels.
\newblock {NLEIGS}: A class of fully rational {Krylov} methods for nonlinear
  eigenvalue problems.
\newblock {\em SIAM Journal on Scientific Computing}, 36(6):A2842--A2864, 2014.

\bibitem{xiao2015contour}
Jinyou Xiao, Shuangshuang Meng, Chuanzeng Zhang, and Changjun Zheng.
\newblock Resolvent sampling based {Rayleigh-Ritz} method for large-scale
  nonlinear eigenvalue problems.
\newblock {\em arXiv preprint arXiv:1510.07522v3}, 2015.

\end{thebibliography}

\end{document}